\documentclass{imsart}

\RequirePackage[OT1]{fontenc}
\usepackage{a4wide}

\usepackage[ruled,linesnumbered]{algorithm2e}
\usepackage{xr-hyper}
\usepackage{hyperref}

\usepackage{amsmath,amsthm,amssymb,mathrsfs}
\usepackage{stmaryrd}
\usepackage{nicefrac}
\usepackage{mathtools} 
\usepackage[utf8]{inputenc}
\usepackage{amsfonts,dsfont}
\usepackage{subfigure}

\usepackage{color}
\usepackage[dvipsnames]{xcolor}
 \usepackage[normalem]{ulem} 

\usepackage{soul,graphicx}
\graphicspath{
  {fig/}
}
\usepackage{natbib}
  \usepackage{pgf, tikz}
  \usetikzlibrary{fit,patterns,decorations.pathreplacing}

\startlocaldefs

\usepackage{amsthm}

\newtheorem{theorem}{Theorem}[section]
\newtheorem{proposition}[theorem]{Proposition}
\newtheorem{corollary}[theorem]{Corollary}

\theoremstyle{remark}
\newtheorem{remark}[theorem]{Remark}
\newtheorem{example}[theorem]{Example}
\theoremstyle{definition}
\newtheorem{definition}[theorem]{Definition}

\usepackage[ruled]{algorithm2e}

\newcommand{\R}{\mathbb{R}}

\renewcommand{\P}{\mathbb{P}}
\newcommand{\E}{\mathbb{E}}

\newcommand{\TDP}{\mbox{TDP}}
\newcommand{\TDR}{\mbox{TDR}}
\newcommand{\FDR}{\mbox{FDR}}
\newcommand{\BH}{\mbox{BH}}
\newcommand{\empBH}{\widehat{\mbox{BH}}}
\newcommand{\empBHsplit}{\widehat{\mbox{BH}}_{\mbox{\tiny Split}}}
\newcommand{\empBY}{\widehat{\mbox{BY}}}
\newcommand{\FDP}{\mbox{FDP}}

\newcommand{\mtc}{\mathcal}

\newcommand{\wt}[1]{{\widetilde{#1}}}
\newcommand{\wh}[1]{{\widehat{#1}}}
\newcommand{\ol}[1]{\overline{#1}}
\newcommand{\ind}[1]{{\mathds{1}_{\{#1\}}}}

\newcommand{\cH}{{\mtc{H}}}

\renewcommand{\l}{\ell}

\newcounter{nbdrafts}
\makeatletter
\newcommand{\checknbdrafts}{
\ifnum \thenbdrafts > 0
\@latex@warning@no@line{**********************************************************************}
\@latex@warning@no@line{* The document contains \thenbdrafts \space draft note(s)}
\@latex@warning@no@line{**********************************************************************}
\fi}

\makeatother

\renewcommand{\subset}{\subseteq}

\endlocaldefs

\usepackage{xr}
\begin{document}

\begin{frontmatter}
%

\title{
{Semi-supervised multiple testing}
}
\runtitle{Semi-supervised multiple testing
}

\begin{aug}
\author{\fnms{David} \snm{Mary}
\ead[label=e1]{david.mary@oca.eu}}
\address{
 Universit\'e C\^ote d'Azur, Observatoire de la C\^ote d'Azur, CNRS, Laboratoire Lagrange, \\ Bd de l'Observatoire, CS 34229, 06304, Nice cedex 4, France \\
\printead{e1}\\
}

\author{\fnms{Etienne} \snm{Roquain}
\ead[label=e2]{etienne.roquain@upmc.fr}}
\address{
Sorbonne Universit\'e (Universit\'e Pierre et Marie Curie), LPSM,\\ 4, Place Jussieu, 75252 Paris cedex 05, France\\
\printead{e2}\\
}
\end{aug}

\begin{abstract}
An important limitation of standard multiple testing procedures is that the null distribution should be known. Here, we consider a 
 null distribution-free approach for multiple testing in the following semi-supervised setting: the user does not know the null distribution, but has at hand a sample drawn from this null distribution. In practical situations, this null training sample (NTS) can come from previous experiments, from a part of the data under test, from specific simulations, or from a sampling process. In this work, we present theoretical results that handle such a framework, with a focus on the false discovery rate (FDR) control and the Benjamini-Hochberg (BH) procedure. 
First, we provide upper and lower bounds for the FDR of the BH procedure based on empirical $p$-values. These bounds match when $\alpha (n+1)/m$ is an integer, where $n$ is the NTS sample size and $m$ is the number of tests. Second, we  give a power analysis for that procedure suggesting that the price to pay for ignoring the null distribution is low when $n$ is sufficiently large in front of $m$; namely $n\gtrsim m/(\max(1,k))$, where $k$ denotes the number of ``detectable'' alternatives. 
Third, to complete the picture, we also present a negative result that evidences an intrinsic transition phase  to the general semi-supervised multiple testing problem {and shows that the empirical BH method is optimal in the sense that its performance boundary follows this transition phase}. 
Our theoretical properties are supported by numerical experiments, which also show that the delineated boundary is of correct order without further tuning any constant.
Finally, we demonstrate that our work provides a theoretical ground for standard practice in astronomical data analysis, and in particular for  the procedure proposed in \cite{Origin2020} for galaxy detection.
\end{abstract}


\end{frontmatter}

\section{Introduction}

\subsection{Background and motivating examples}\label{sec:background}

Multiple testing, with emphasis on large scale problems, is an important topic in modern statistics. Classical theory and performance guarantees heavily rely on the knowledge of the null distribution. However, in many practical situations, the null distribution is out of reach. 
A famous situation, described in a series of work by \cite{Efron2004, Efron2007b,  Efron2008,Efron2009b} and followed by, e.g., \cite{Sch2010, AS2015, Ste2017,SS2018,roquain2020false} is the case where the null distribution is mis-specified and is empirically adjusted from the data by fitting some parametric null model (typically Gaussian). In particular, it is well known that using an erroneous null can by disastrous in terms of false discovery rate (FDR), see, e.g.,  \cite{RVvignette2020}. 
Related works, relying on the famous two-group model (\citealp{ETST2001}), propose to estimate the null distribution together with the proportion of nulls and the alternative distribution, and to plug them into the so-called local FDR values, see \cite{ETST2001} and \cite{PB2012,heller2014} among others. The latter can in turn be used into an FDR controlling procedure, see \cite{SC2007,SC2009,CS2009,
CSWW2019,roquain2020false,abraham2021multiple}. The validity of such approaches, often given asymptotically in the number of tests, also requires strong model assumptions to ensure that these parameters  can be correctly estimated. 

Here, we consider a semi-supervised setting, with essentially no assumption on the null distribution. Instead, the user has at hand a sample, called the {\it null training sample} {(NTS)}, of length $n\geq 1$, and generated according to this unknown null.  This is motivated by the two following generic situations:
\begin{itemize}
\item {\it Blackbox null sampling:}  the exact expression of the null distribution is intractable, but a sampling machine is able to simulate according to the null distribution. In that case, the NTS is exogenous and its length $n$ corresponds to the number of sampling, so can be chosen by the user. It is nevertheless typically limited in size by computation time constraints. 
\item {{\it Null sample given:}  the null distribution is unknown, but previous experiments or experts provide a fixed number $n$ of examples under the null. The NTS is exogenous as in the above case, but $n$ cannot be modified by the user.  }
\item {{\it Null sample learned from data:}} the null distribution is unknown, but an independent part of the same data set provides an NTS for the user.
In that case, the NTS is endogenous, of a given length $n$ that cannot be modified by the user. 
\end{itemize}

The case of ``blackbox null sampling'' is motivated by numerous  situations. Two motivations come from Astrophysics; first when a code can be used to simulate images of astrophysical  sources, see e.g. \cite{Cosmic2021} (their Figure 15). Second, when the NTS comes from instrumental captures that are made without the objects of interest, see e.g.  \cite{Choquet_2018} for the detection of exoplanetary debris disks (their Figure 5). In each of these situations, the null distribution is not accessible for the user, and only the NTS can be generated. More broadly, this case is motivated by recent advances in machine learning, especially implicit generative models, as generative adversarial networks (\citealp{Good2014}), or variational auto-encoders (\citealp{Kingma2014}), for which sampling is possible without knowing the underlying distribution. An illustration of the blackbox null sampling case is provided in Section~\ref{sec:MC}, on a toy example for which  multiple likelihood ratio tests are simultaneously performed. 

{The case of ``null sample given''  is common in the machine learning context, where the learner is given a sample of ``nominal patterns'' but without labeled novelties. This is classically referred to as ``one class classification'' or 
``learning from positive and unlabeled examples'' and we refer the reader to the work \cite{BLS2010} that pointed out many references in this abundant literature.
}

The case of ``null sample learned from data'' refers to the framework where it is possible to isolate part of the data to produce a sample that contains copies of the test statistics under the null, or approximately so. 
While it can be met in various datasets, it is motivated by a specific application in Astrophysics that is extensively developed in Section~\ref{sec:appli}. 
It regards the detection of galaxies in the early Universe from image measurements in multiple wavelength channels. In this application, the distribution of the tests statistics under the null is  unknown and it was proposed in \cite{Origin2020, Cosmic2021} to estimate this distribution from a null training sample obtained from the data itself. The NTS is obtained as the population of the opposite of local minima and the whole NTS is used for testing each of the $m$ local maxima.

In both cases, a crucial issue is to build a procedure for making discoveries while being fully interpretable, especially when the number of tests $m$ is large. We thus focus on building a procedure that controls the false discovery rate (FDR), that is, the expected ratio of errors among the discoveries made by the procedure (\citealp{BH1995}). 
Interestingly, controlling the FDR by using a simulated NTS has similarities 
with the recent ``knockoff'' method introduced in \cite{BC2015} which has been at the origin of an impressive scientific production over the last years, see, e.g., \cite{weinstein2017power,katsevich2019multilayer,barber2019knockoff,bates2020metropolized}. Further comparisons are given in Section~\ref{sec:relatedwork}.

When proper $p$-values can be built, the classical way to control the FDR at level $\alpha$ is to use the Benjamini Hochberg (BH) procedure (\citealp{BH1995}). However, in the setting described above, the exact $p$-values are out of reach, so that the usual BH procedure cannot be used. In our context, we call it the {\it oracle} BH procedure, and denote it by $\BH^*_\alpha$, or $\BH^*$ for short.  
Instead, the NTS can be used to build empirical $p$-values, called $\hat{p}$-values for short.
It is then natural to use the $\hat{p}$-values into the BH procedure, which is the procedure studied in this paper. We call it the {\it semi-supervised} BH procedure and denote it by $\empBH_\alpha$, or $\empBH$ for short.

Let us  already note that plugging empirically-based $p$-values into the BH procedure is not new and has been widely explored in the literature, especially in a Monte Carlo framework, see, e.g.,  \cite{GP2008,sandve2011sequential,GH2014,zhang_adaptive_2019}. However, while the same null sample is used to compute all $p$-values in our setting, most of the existing works focus on the case where $m$ null samples are available, that is, each test uses a different sample, often generated via randomization process (e.g., permutations). In that case, the computational price is much higher and these works mostly aim at reducing this price. 
The case of {\it only one} null sample has been considered only recently to our knowledge, see \cite{weinstein2017power,bates2021testing}. 
The computational issue can be easily solved (see Algorithm~\ref{algo:empBH}), and our emphasis is rather on the theoretical guarantees of the resulting BH procedure ($\empBH$). Further details and comparisons with existing literature are given in Section~\ref{sec:relatedwork} and in Section~\ref{sec:naive} {in the supplement}.

 Finally, an important point of our work will be to determine how large  $n$ should be relatively to the number $m$ of tests.  
 Obviously, when $n$ tends to infinity while the number $m$ of tests is kept fixed, the situation becomes similar to the one where the null distribution is known (that is, when $n=\infty$). 
But the situation is more complex when both $n$ and $m$ gets large simultaneously, which is {typical} (e.g., in our galaxy detection example, we have $n\approx m = 3.3\times 10^6$). 
As can be guessed, the full picture also depends on the sparsity of the signal. This will be adressed in our theory through a parameter called $k$, which is a proxy to the number of detectable alternatives.

\subsection{Contributions}

\paragraph{Main contributions.} The main contributions of the paper can be summarized as follows:

First,  we study the FDR of the procedure  $\empBH$, by providing upper and lower bounds (Theorem~\ref{thFDR}). 
These bounds hold in a strong sense, that is, for any couple $(n,m)$ with $n,m \geq 1$, any number of true nulls $m_0$, any null distribution, and any marginal distribution of the alternatives.
Moreover,  these bounds match and equal $\alpha m_0/m$ when $\alpha (n+1)/m$ is an integer.  In practice, this provides a first guideline for choosing $n$ in order to avoid over-conservativeness of the procedure.

Second, 
we provide a power boundary 
 for $\empBH$, which puts forward the crucial role of $n$ with respect to $m$:   the power of $\empBH$ is close to the one of the oracle $\BH^*$ if $n\gtrsim m/\alpha$ (Proposition~\ref{prop:power}), but is not when $n\lesssim m/\alpha$ (Proposition~\ref{thm:lb0}). This leads to the boundary $n\asymp m/\alpha$.
In addition, we underline the role of the sparsity in the boundary with the following additional result. {For distributions that are more favorable in the sense that the oracle $\BH^*$
 is expected to make at least $k$ true discoveries with high probability (a situation where we say that $k$ alternatives are ``detectable'')}, we show that the power boundary for $\empBH$ occurs at 
   $n\asymp m/(k\alpha)$. 
 As an illustration, for $k=1$, the boundary is $n\asymp m/\alpha$ and thus is the same as for general distributions. However, for the dense case $k=m/2$, the boundary reads $n\asymp 1/\alpha$. This indicates that an NTS of size $\gtrsim 1/\alpha$ is enough to recover the power of the oracle in this case. {This is markedly different from the case of general distributions. In particular, oracle performances can be achieved in the dense case for a constant value of $n$, regardless of $m$.}
Overall, this leads to a new ``rule of thumb'' with a transition at $n=m/(\alpha\max(1,k))$, which is implemented in the numerical experiments, see Section~\ref{sec:num} and in the astrophysical example, see Section~\ref{sec:appli}.

Third,  we show that an intrinsic phase transition occurs in the general case at $n\asymp m$ (Corollary~\ref{cor-lb}). The boundary $n\asymp m$ ($\alpha$ being fixed) can not be improved by another procedure: when $n\lesssim m$, 
 {no procedure (only based on the observations and the NTS)} can both control the FDR while having a power close to the one of $\BH^*$   (Theorem~\ref{thm:lb}). 
Since $\empBH$ does mimic the oracle when $n\gtrsim m$ (Proposition~\ref{prop:power}), this 
establishes a general minimax-type optimality property for $\empBH$.
(Note that the test statistic is fixed in our setting so that $\BH^*$ is an appropriate reference for power, see Section~\ref{sec:BHoptimal} for a further discussion.)


\paragraph{Secondary contributions. }
{Additional secondary contributions are as follows:}


First, we show how $\empBH$ can be used in the ``Blackbox null sampling'' setting in Appendix~\ref{sec:MC}. 
We introduce the Blackbox BH procedure, which is defined as the semi-supervised BH procedure with a preliminary step where the NTS is properly generated, see Algorithm~\ref{algo:MCBH}.
While it can be used in a very broad context, we illustrate its use for likelihood ratio tests for which the oracle is accessible in Section~\ref{sec:LRT} {in the supplement}. A comparison with local FDR type approaches is also provided in that case.

Second, 
we put forward the following, perhaps seemingly paradoxal, fact for FDR control under negative dependence. Even in the classical setting where the true null is known, it is better not to use BH procedure, but to build instead  artificially an NTS, and to use it along with the semi-supervised procedure $\empBH$. This approach is refered to the randomized BH procedure, which is studied separately in Appendix~\ref{sec:randBH}. While the superiority of the randomized BH procedure over the usual BH procedure in terms of FDR control is shown 
for an admittedly restrictive dependence structure, correcting the BH procedure to accommodate negative dependencies is known to be a challenging task (see, e.g., \cite{fithian2020conditional} and references therein). We think that this intriguing side result is an important proof of concept for the randomized BH procedure.

Third, extensive numerical experiments are given in Section~\ref{sec:num} that validate and illustrate our theoretical results. In particular, they corroborate the fact that the boundary where the power of $\empBH_\alpha$ gets of the order of the one of $\BH^*_\alpha$ occurs around $n= m /(k\alpha)$ (without further tuning of the constant), where $k$ is the number of ``detectable'' alternatives in the data.  {For instance, and perhaps counter-intuitively, it is shown that oracle performances can be achieved in a dense case {for values of $n$ as small as $5$ or $10$, regardless of $m$.}

Fourth, a detailed application to galaxy detection is given in Section~\ref{sec:appli}.
{Remarkably, the recent results of \cite{Cosmic2021} suggest
the likely discovery of an unexpected population of ultra-faint dwarf galaxies\footnote{{
Also disseminated by the CNRS press release, see,}  \href{https://www.cnrs.fr/en/first-images-cosmic-web-reveal-myriad-unsuspected-dwarf-galaxies}{https://www.cnrs.fr/en/first-images-cosmic-web-reveal-myriad-unsuspected-dwarf-galaxies}}. This discovery results  from a two-stage detection process, whose first stage relies on {a {former} version of the semi-supervised Benjamini-Hochberg procedure developed in  \cite{Origin2020}, which also provides the same output as $\empBH_\alpha$. Hence, the } present paper provides a theoretical support to {these findings}, with guarantees both on the FDR and on the power. }\\

Figure~\ref{fig:phasetransition} summarizes the different power regimes put forward in our analysis. The transition phase $n=m/\alpha$ separates two regimes: the regime where oracle performances can be reached for any distribution (``mimicking the oracle possible in general'', lime green) versus the regime where no procedure can reach the oracle performances (``mimicking the oracle impossible in general'', tomato $+$ red brick).  The line $n=m/(\alpha k)$ is the performance boundary of $\empBH$ for favorable distributions for which at least $k$ alternatives are detectable (in which case oracle performances can be reached in the lime green $+$ tomato area). 
Note that our theory proves that these boundaries hold only up to numerical constants, whereas the numerical experiments suggest that they hold with constant $1$.

\begin{figure}[h!]
    \centering
\includegraphics[width=.5\textwidth]{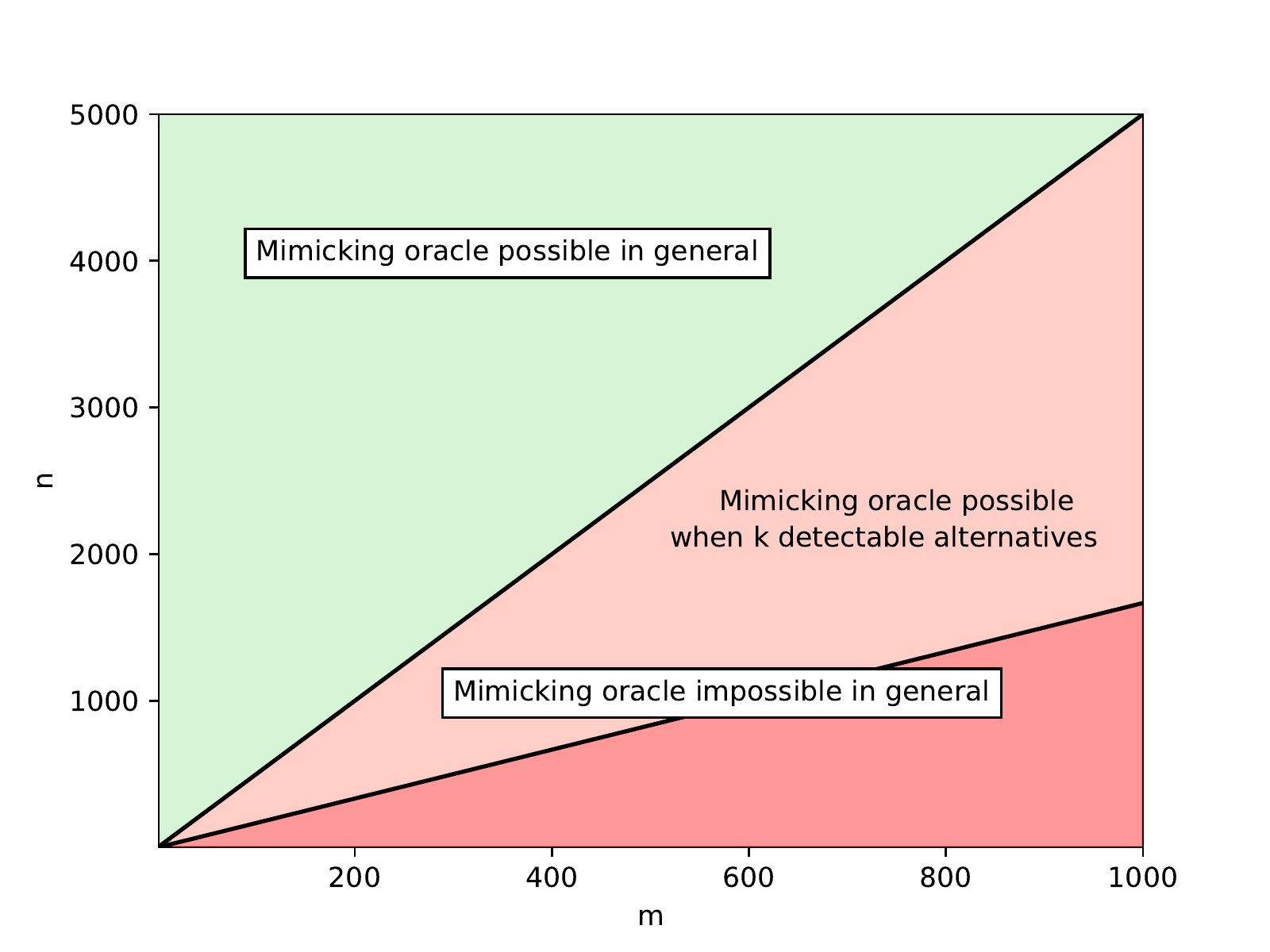}  \hspace{-5mm}\includegraphics[width=.5\textwidth]{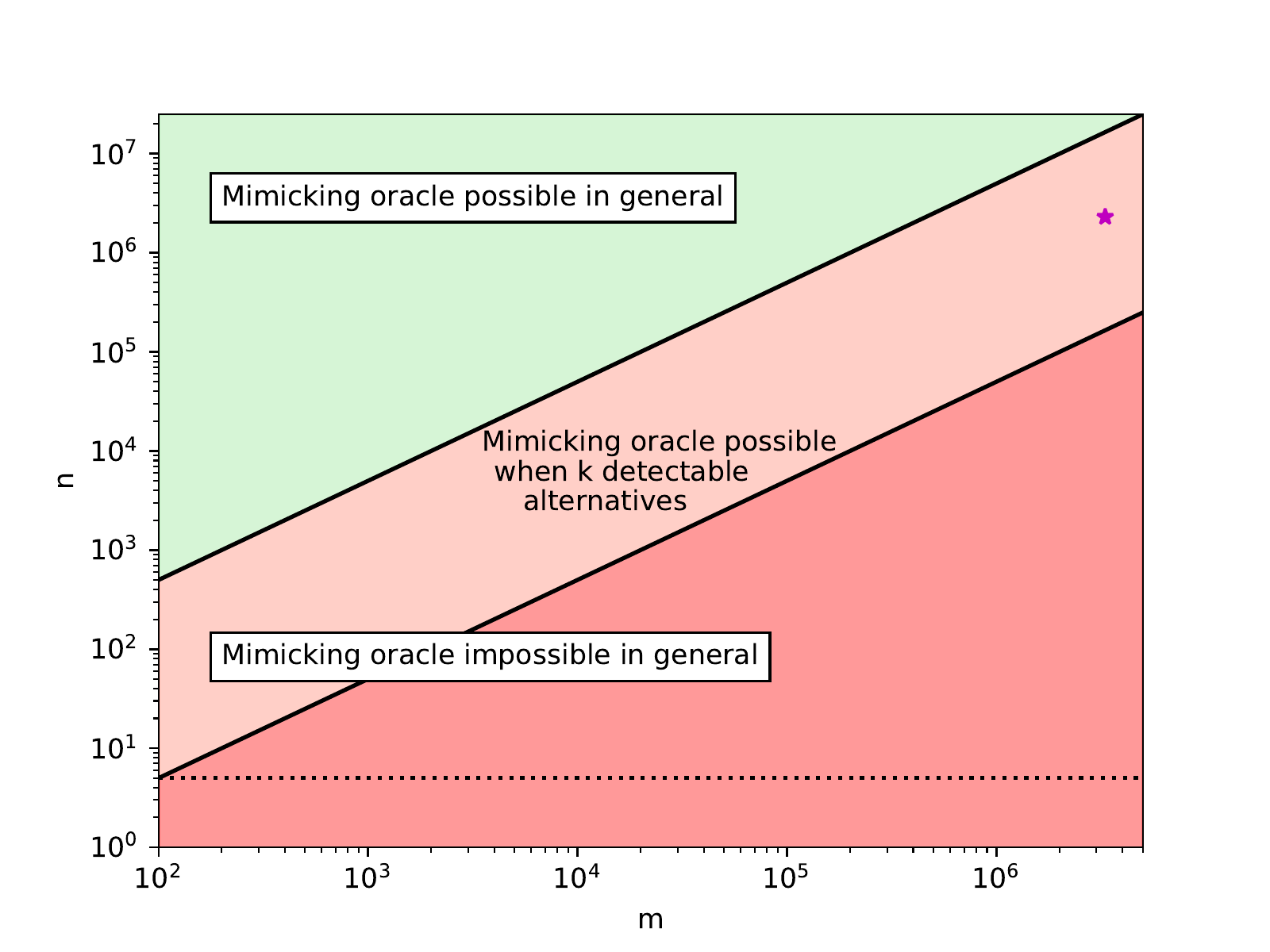} 
      \caption{\label{fig:phasetransition}  Visualization of the general, distribution-free phase transition $n=m/\alpha$ for the semi-supervised multiple testing problem, as established in Section~\ref{sec:boundary} (with $\alpha=0.2$) and of the $\empBH$ boundary $n=m/(\alpha k)$, only valid for distributions with at least $k$ detectable alternatives  in the sense defined in Section~\ref{sec:morefavorable}. Left:   $k=3$. Right: $k=100$; plot in $\log$-$\log$ scale; the  boundary $n=1/\alpha$ ($k=m$) is added with a dotted line;  the case of the MUSE data set (Section~\ref{sec:appli}) is also added with a star symbol.} 
        \end{figure}

\subsection{Related works}\label{sec:relatedwork}

\paragraph{Permutation-based multiple testing.}

A common way to generate a ``null sample'' from the data under test is to apply some randomization that preserves the null distribution, typically by performing permutations of individuals. 
While single testing using randomization is  classical and can be traced back to \cite{Fish1935}, several extensions have been proposed in the literature to accommodate multiple testing criteria, see \cite{WY1993,lin2005efficient,RW2005,RW2007,HSG2019}.
In particular, an active line of research is dedicated to reduce the computation time of BH procedure with ${p}$-values obtained from permutation-based null samples: indeed, the usual permutation-based paradigm requires to generate a different ``null sample'' for each test, which makes the use of such a BH procedure prohibitive in that framework. 
In \cite{GP2008}, they adapt the number of bootstrap samples sequentially to speed-up BH procedure by using 
bootstrap confidence intervals for $p$-values. This method is further refined in \cite{GH2014}, where  the procedure recovers with high probability the rejection set of the BH procedure using  ``ideal'' $p$-values (exhausting all permutations). 
Another approach is used in \cite{sandve2011sequential}  by allocating the Monte Carlo budget (total number of Monte Carlo samples) according to the significance of the test statistics, itself extending an idea of \cite{besag1991sequential} for single testing. More recently, \cite{zhang_adaptive_2019} proposed to reduce the computation burden by following a bandit approach.
While all these works are based on null training samples, the crucial difference is that our setting only relies on  {\it one} null sample for all tests. The consequences are the following: first, the complexity of the procedure proposed here ($\empBH$) is much smaller than that of the  BH procedure with permutation-based $p$-values, 
the need for designing an efficient algorithmic strategy is far less critical than
 in the  works mentioned above (note that our Algorithm 1 for $\empBH$ is nevertheless efficient). Second, this advantage comes with a counterpart: in the case where the initial test statistics are independent, the permutation-based ${p}$-values are also independent, 
 while our setting induces dependencies between the $\hat{p}$-values (the same NTS is used to build all $\hat{p}$-values). This makes the FDR control more difficult to obtain. 

 Finally, the above comparison has to be moderated by the fact that randomization testing and our semi-supervised setting each come with specific mathematical assumptions:  randomization testing relies on a null distributional invariance which is very different from the assumption \eqref{exch} below. Namely, the exchangeability property concerns the set of ``variables" (nulls of the test sample  plus  the null training sample), whereas in permutation testing, the exchangeability concerns the set of individuals.  
As a result, mathematical results derived in each framework cannot be directly compared. In particular, it is important to note that we do not pretend to address the FDR controlling problem in the permutation-based framework. Our contribution lies in another framework, which thus departs from the  Monte-Carlo literature mentioned above.

\paragraph{Earlier occurrences of $\empBH$.}

 In the linear Gaussian model, \cite{BC2015} proposed to build test statistics (from so-called ``knockoff'' variables) that have a special symmetry property under the null allowing to properly calibrate an FDR controlling procedure.
Still in that framework, a rapidly growing literature proposed further extensions and refinements of this seminal work, see, e.g.,
\cite{candes2018panning,barber2019knockoff,bates2020metropolized,Barber2020robust,liu2018auto,nguyen2020aggregation,sarkar2021adjusting}. 
Here, the empirical BH procedure can be seen as an extension of \cite{BC2015} to the semi-supervised setting. 
It turns out that this procedure has been considered in the paper by \cite{weinstein2017power}, whose scope is yet quite different. Very recently, it has been considered by \cite{bates2021testing,yang2021bonus}, for which the NTS is called the inlier sample and the bag of nulls, respectively. These papers come with theoretical guarantees, see Remark~\ref{rempreviousFDR} below. Our results have been developed independently.

\paragraph{Multiple comparisons to control.}

{
Multiple comparisons to control (MCC) is a long-established problem in multiple testing \citep{dunnett1955multiple,hsu1996multiple,FS2007b,fithian2020conditional} where one typically aims at comparing several treatments to some common benchmark (control). In the MCC setting, one typical observes only one test statistic per treatments and one test statistics for the control. This would correspond to the case where the null training sample is of length $n=1$, which is not the typical case considered here. 
Hence, to our knowledge, the connection to that part of the literature is only weak.
}

\paragraph{Other FDR controls.}

Our work is closely related to the task of semi-supervised novelty detections \citep{BLS2010}, developed in a machine learning context, where the user has at hand both a null sample and an unlabeled sample and they aim at labeling the unlabeled sample. However, the procedures developed therein are significantly different from here: first, they adjust the test statistics by considering families of classifiers. 
Second, their FDR control is based on a concentration argument that adds an error term larger than $n^{-1/2}+m^{-1/2}$ (see Proposition 12 therein) and depending on the VC-dimension of the classifier class, while the FDR control in Theorem~\ref{thFDR} is exact (no error term). 
 
 Finally, another closely related literature tackles the issue of learning the null distribution without null training sample (only using the original test statistics) but assuming that the null distribution belongs to a parametric model, typically Gaussian with unknown mean and variance. While the most classical line of research is the one following the ``local FDR" methodology introduced by Efron, see, e.g., \cite{Efron2008}, 
 theoretical results have been obtained by  \cite{CDRV2021,roquain2020false}.  
 The methodology developed here, and particularly the impossibility result (Section~\ref{sec:impos}) and the boundary phenomenon (Section~\ref{sec:boundary}), are inspired from \cite{roquain2020false}. However, the setting being markedly different, several substantial adjustments are required. Also, we underline that we derive here an FDR control without remainder terms, which was not the case in  \cite{CDRV2021,roquain2020false}.

\paragraph{Naive solutions to our problem.}
For completeness, let us discuss two naive solutions that can be straightforwardly used to derive a procedure with a proven FDR control in the present semi-supervised setting, and explain why they are not satisfactory. Recall that, even under independence of the test statistics, the $\hat{p}$-values are not independent, which is a problem to design an FDR controlling procedure that takes as input these $\hat{p}$-values.

First, one solution is to use the Benjamini-Yekutieli  procedure or one of its extension \cite{BY2001,BR2008} that control the FDR under arbitrary dependence between the $p$-values, so also when used with $\hat{p}$-values. Namely, the semi-supervised Benjamini-Yekutieli procedure, denoted by $\empBY_\alpha$ (or $\empBY$ for short),  considers $\empBH_{\alpha/c_m}$ at level $\alpha/c_m$ where $c_m=1+1/2+\dots+1/m$. However, it is well known that the power loss is substantial with respect to BH procedure 
{and this general fact also holds  in our setting, as it will be shown in the numerical experiments, see Section~\ref{sec:naive} {in the supplement}}.
In addition, Theorem~\ref{thFDR} shows that under Assumption~\eqref{exch}, the procedure $\empBH$ already achieves the desired FDR control so there is no need to use the corrected procedure $\empBY$.

A second naive solution, referred to as $\empBHsplit$, is to split the NTS of size $n$ into $m$ null samples $T^1,\dots,T^{m}$, each of size $n/m$ (say that the latter ratio is an integer for simplicity) so that each $\hat{p}$-value uses a different part of the null sample, that is, each $\hat{p}_i$ is computed from the null training sample $T^i$.
In that case, if the test statistics are independent, these modified  $\hat{p}$-values are also  independent, and the BH procedure using these modified $\hat{p}$-values does control the FDR by the original result of \cite{BH1995}. However, this reduces drastically the size of the (different) NTS used to calibrate each test ($n/m$ instead of $n$), which leads again to a poor power, see Section~\ref{sec:naive}. 

\subsection{Organization of the paper}

The paper is organized as follows: while the model, procedures and criteria are detailed in Section~\ref{sec:prelim}, the FDR results are given in Section~\ref{sec:FDR}.
Power properties of $\empBH$ are then derived in Section~\ref{powerresult} with upper and lower bounds, which delineate boundaries for $\empBH$. Extending to any procedure the impossibility result below the boundary, the result of Section~\ref{sec:optimality} delivers an optimality property of $\empBH$ and a general phase transition for the semi-supervised multiple testing problem.
We then illustrate our findings with numerical experiments in Section~\ref{sec:num} and the motivating application to astrophysical data is investigated in Section~\ref{sec:appli}. We conclude and discuss several open issues related to our work in Section~\ref{sec:conclusion}.
Two by-products of our theory are presented in Appendices~\ref{sec:MC}~and~\ref{sec:randBH}, with the blackbox BH procedure and the randomized BH procedure, respectively. 
{For space reasons, we have deferred some materials to a supplemental file, whose sections are numbered with the prefix ``S'' to avoid any confusion. In this supplement, the main proofs are given in Section~\ref{sec:proofs} and Section~\ref{sec:proofspower} for the FDR results and the power results, respectively.  
Auxiliary results and proofs are postponed to Section~\ref{sec:auxiliary}, while additional numerical experiments are given in Section~\ref{add:num}.}

\section{Preliminaries}\label{sec:prelim}

\subsection{Setting}\label{sec:setting}

For $n,m\geq 1$, let us observe a sample
$Z=(Z_1,\dots,Z_{n+m})=(Y_1,\dots,Y_n,X_1,\dots,X_m)\in \R^{n+m},$ whose distribution is denoted by $P$, the model parameter, that belongs to some model $ \mathcal{P}$. 
The sample $Y=(Y_1,\dots,Y_n)$ is referred to as the null training sample (NTS), which is assumed to be identically distributed of marginal distribution $P_0=P_0(P)$. We denote the upper-tail function of $P_0$ by $F_0(t)=\P(X_i\geq t)$, $t\in\R$, $i\in \cH_0(P)$, which is assumed to be continuous and decreasing on the support of $P_0$. This will be the only assumption made on $P_0$  throughout the manuscript.

The sample $X=(X_1,\dots,X_m)$ corresponds to the sample under test. 
We consider the multiple testing problem where we would like to test the $i$-th null hypothesis $H_i$: ``$X_i\sim P_0$" (against the complementary alternative), simultaneously for $1\leq i\leq m$. {Note that while we allow for arbitrary alternatives here, this setting is typically suitable for alternatives that make $X_i$ stochastically larger than under the null (decisions will be based upon large values of the $X_i$'s).} 
Classically, let us denote $\cH_0(P)=\{i\in\{1,\dots,m\}\::\: X_i\sim P_0\}\subset \{1,\dots,m\}$ the subset corresponding to true null hypotheses and $m_0(P)=|\cH_0(P)|$.
Let us denote $\cH_1(P)$ the complement of $\cH_0(P)$ in $\{1,\dots,m\}$ and  $m_1(P)=m-m_0(P)$. Often, we omit the parameter $P$ in the notation $P_0,\cH_0,\cH_1,m_0,m_1$ for simplicity.

Throughout the paper, we are going to consider various dependence assumptions between the $Z_i$'s. The most simple assumption is
\begin{align}\label{indep}\tag{Indep}
\mbox{$(Y_1,\dots,Y_n,X_i,i\in\cH_0)$ are i.i.d. $\sim P_0$ and independent of  $(X_i,i\notin\cH_0)$}.
\end{align}
Note that \eqref{indep} does not exclude dependencies between the elements of $(X_i,i\notin\cH_0)$.
We also use the following less restrictive condition:
\begin{align}\label{exch}\tag{Exch}
\mbox{$(Y_1,\dots,Y_n,X_i,i\in\cH_0)$ are exchangeable and independent of $(X_i,i\notin\cH_0)$}.
\end{align}
Hence, under \eqref{exch}, there could be also some dependencies  between the elements of $(Y_1,\dots,Y_n,X_i,i\in\cH_0)$.

\subsection{Procedures, criteria and $p$-values}
\label{sec:pcp}
A multiple testing procedure is a (measurable) function $R=R(Z)$ that returns a subset of $\{1,\dots,m\}$ corresponding to the indices $i$ where $H_i$ is rejected. For any such  procedure $R$, the false discovery rate (FDR) of $R$ is defined as the average of the false discovery proportion (FDP) of $R$ under the model parameter $P\in \mathcal{P}$, that is, 
\begin{align}
\FDR(P,R)&= \E[\FDP(P,R)],\:\:\: \FDP(P,R)=\frac{\sum_{i\in \cH_0} \ind{i\in R}}{1\vee |R|}.\label{equFDRFDP}
\end{align}
{Similarly, the true discovery rate (TDR) is defined as the average of  the true discovery proportion (TDP), that is,}
\begin{align}
\TDR(P,R)&= \E[\TDP(P,R)],\:\:\: 
\TDP(P,R)=\frac{\sum_{i\in \cH_1} \ind{i\in R}}{1\vee m_1(P)}.\label{equTDRTDP}
\end{align}
Note that if $m_1(P)=0$, $\TDP(P,R)=0$ for all procedures $R$.

In the sequel, we will focus on $p$-value based procedures and we implicitly consider the situation where it is desirable to reject $H_i$ for large values of $X_i$.
 If the null distribution $P_0$ is known,  $F_0$ is known and we can consider $p_i(X)=F_0(X_i)$, $1\leq i\leq m$. 
By definition, the $p$-value family $p_i=p_i(X)$, $1\leq i\leq m$, satisfies that for all $ i \in \cH_0(P)$, $p_i\sim U(0,1)$, and thus also the super-uniformity property 
\begin{equation}\label{equ:superunif}
\forall i \in \cH_0(P), \:\:\forall u\in [0,1], \P_{Z\sim P}(p_i\leq u)\leq u.
\end{equation}
As it is required to obtain valid individual tests, condition \eqref{equ:superunif} is generally considered as the definition of ``valid'' $p$-values.

Since in our framework $P_0$ is unknown, the above $p$-values are unknown oracle $p$-values and thus cannot be used in practice. Instead, the null sample $(Y_1,\dots,Y_n)$ can be used to build the empirical  $p$-values 
\begin{equation}\label{equptilde}
\wt{p}_i(Z) = n^{-1}\sum_{j=1}^n \ind{Y_j\geq X_i},  \:\: 1\leq i\leq m.
\end{equation}
However,  the $\wt{p}_i $'s do not satisfies the necessary super-uniformity \eqref{equ:superunif}. For instance, for $u=0$, the condition \eqref{equ:superunif} is violated because the {event} $\wt{p}_i(Z)=0$ can occur with a positive probability.
{Hence, using the $\wt{p}_i $'s as $\hat{p}$-values is not appropriate, especially in a multiple testing context where under-estimating $p$-values can lead to an increased number of false discoveries. This phenomenon is well known and we refer the reader to the review of \cite{phipson2010permutation} for more details on this issue (see also the references therein). }
 A common way to correct the $\wt{p}_i $'s is to make them slightly biased upward  by considering rather the conservative version (see, e.g.,  \citealp{davison1997bootstrap}), given by
\begin{equation}\label{equemppvalue}
\wh{p}_i(Z) = \wh{F}_0(X_i) =(n+1)^{-1}\sum_{x\in \{X_i,Y_1,\dots,Y_n\}} \ind{x\geq X_i}, \:\: 1\leq i\leq m,
\end{equation}
where we let
\begin{equation}\label{equF0chap}
\wh{F}_0(x) = (n+1)^{-1}\left(1+\sum_{j=1}^n \ind{Y_j\geq x}\right), \:\:\:x\in\R.
\end{equation}
Under \eqref{exch}, since for any $i\in \cH_0$, the variables $X_i,Y_1,\dots,Y_n$ are exchangeable, the $\wh{p}_i(Z)$'s do satisfy the super-uniformity \eqref{equ:superunif}, see, e.g., Lemma~5.2 in \cite{ABR2010a}. 
Hence, the $\wh{p}_i(Z)$ are ``valid'' $p$-values, that can in turn be plugged into multiple testing procedures.

\subsection{BH procedures}

In this work, an important class of multiple testing procedures is the BH-type procedures, which use as input different $p$-value families.  
The BH procedure  
is defined as follows: 
for some level $\alpha\in(0,1)$, order the $p$-values in increasing order $ p_{(1)}\leq \dots\leq p_{(m)}$ and then let
\begin{equation}\label{equBH}
\BH_\alpha =\{ i\in \{1,\dots,m\}\::\: p_i\leq \alpha \wh{k}/m\} , \:\:\wh{k}=\max\{k\in\{0,1,\dots,m\}\::\:p_{(k)}\leq \alpha {k}/m\},
\end{equation}
where $\alpha$ is the nominal level of BH procedure and where we let $p_{(0)}=0$ by convention.

\begin{definition}\label{def:BHs}
We consider the two following versions of BH procedure, depending on which $p$-value family is given as input:
\begin{itemize}
\item the {\it oracle BH procedure}, denoted by $\BH^*_\alpha$, is the BH procedure using the unknown $p$-values $p_i(X)=F_0(X_i)$, $1\leq i \leq m$;
\item the {\it semi-supervised BH procedure}, denoted by $\empBH_\alpha$, is the BH procedure using the $\hat{p}$-values $\wh{p}_i(Z)$, $1\leq i \leq m$, given by \eqref{equemppvalue}.
\end{itemize}
\end{definition}

\begin{algorithm}
\KwData{$Z=(Z_1,\dots,Z_{n+m})=(Y_1,\dots,Y_n,X_1,\dots,X_m)\in \R^{n+m}$ semi-supervised sample, $\alpha$ level}
 \begin{enumerate}
 \item Order the $Z_i$'s, that is, 
$
Z_{\tau(1)}\geq \dots \geq Z_{\tau(n+m)}
$, for some permutation $\tau$ of $\{1,\dots,n+m\}$
\item Let $s_\l=\ind{\tau(\l)\leq n}\in\{0,1\}$ which is $1$ if and only if $Z_{\tau(\l)}$ comes from sample $Y=(Y_1,\dots,Y_n)$
\item Let $\FDP=1$, $V=n$, $\l=m+n$, $K=m$
\item While ($\FDP>\alpha$ and $K\geq 1$) do $\l=\l-1$ 
\begin{itemize}
\item if $s_{\l+1}=1$, $V=V-1$
\item else, $K=K-1$
\end{itemize}
$\FDP=\frac{{V+1}}{n+1}\frac{m}{K}$ (or $\FDP=1$ if $K=0$)
\end{enumerate}
\KwResult{$\empBH_\alpha=\{i\in \{1,\dots,m\}\::\: X_i\geq X_{(K)} \}$ (reject nothing if $K=0$).}
\caption{$\empBH_\alpha$, the semi-supervised BH procedure}\label{algo:empBH}
\end{algorithm}

\begin{figure}[h!]
\begin{center}
\includegraphics[scale=0.55]{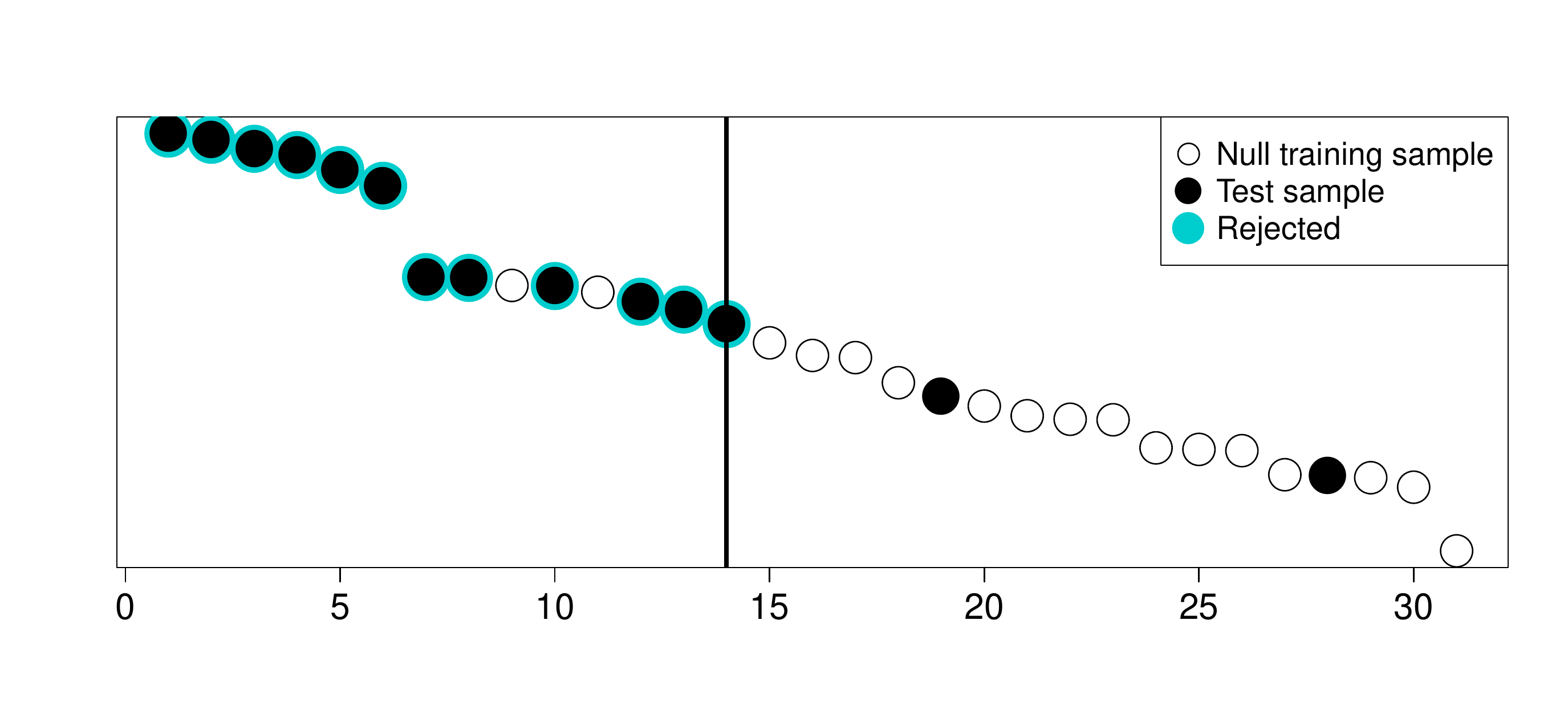}
\end{center}
 \vspace{-1cm}
 \caption{\label{fig:algo} Illustration of Algorithm~\ref{algo:empBH} for $m=14$, $n=17$, $\alpha=0.2$, and some realization of the ordered test statistics. At the point of rejections (vertical line) $\l=14$, we have $V=2$, $K=12$, $\FDP=\frac{2+1}{n+1}\frac{m}{12}\leq 0.2$, while $\FDP>0.2$ in any further point $\l>14$. The algorithm makes $K=12$ rejections in the test sample (depicted in blue). See text for further comments. 
 }
 \end{figure}
 
 
Importantly, the output of $\empBH_\alpha$ can be quickly derived, even for large values of $n$ and $m$, thanks to an algorithm of complexity $O((n+m)\log(n+m))$, see Algorithm~\ref{algo:empBH}. This comes from the reformulation of $\empBH_\alpha$ given in Section~\ref{sec:reform}, which was also used by \cite{weinstein2017power} 
 in another context. Figure~\ref{fig:algo} provides an illustration of Algorithm~\ref{algo:empBH}: it is a stepwise procedure that goes from the smallest values of the test statistics (right) to the largest values (left), and that stops the first time where the FDP falls 
below $\alpha$. At each step, the FDP is estimated by the ratio of the number of null samples in the left part plus one ($V+1$), to the number of test statistics in the left part ($K$), this ratio being sample-sized corrected by the factor $m/(n+1)$. 
Hence, at each step, the $Y_i$'s are used as benchmarks to evaluate how many false discoveries are expected among the considered $X_i$'s.
{Finally, while the above version of Algorithm~\ref{algo:empBH} was presented for simplicity, a shortcut (faster) version can obviously be obtained 
by iterating in the loop only over the indices $\l$ corresponding to the $X_i$'s  (the FDP is computed only at black points in Figure~\ref{fig:algo}).
}

\section{FDR control}\label{sec:FDR}

In this section, we study the FDR of the procedure $\empBH$.

\begin{theorem}\label{thFDR}
For all $n,m\geq 1$ and $\alpha\in (0,1)$, consider the semi-supervised BH procedure $\empBH_\alpha$ at level $\alpha$ as defined in Definition~\ref{def:BHs}. 
Then, for any parameter $P$ satisfying \eqref{exch}, the following holds:
$$
\frac{m_0}{m}\frac{m}{ n+1}\bigg\lfloor \alpha \frac{n+1}{m} \bigg\rfloor \leq \FDR(P,\empBH_\alpha)\leq \alpha m_0/m,
$$
where $\lfloor x \rfloor$ denotes the largest integer smaller than or equal to $x$.
In particular, when $\alpha (n+1)/m$ is an integer, the FDR bound is achieved, that is, $\FDR(P,\empBH_\alpha)= \alpha m_0/m$.
\end{theorem}

The proof is given in Section~\ref{sec:proofs} and is based on a super-martingale argument which is similar to that of \cite{BC2015}. However, a major difference is that the underlying process is not an i.i.d. Bernoulli process, but is only exchangeable, see Lemma~\ref{lem:annexe} for more details. 
The lower bound part is obtained by looking carefully at the remainder term in the super-martingale property. To our knowledge, this kind of refinement is new in the  literature. This allows to evaluate the sharpness of the FDR bound.


{In particular, }Theorem~\ref{thFDR} shows that under \eqref{indep} (implying \eqref{exch}) the semi-supervised BH procedure $\empBH_\alpha$ has an FDR smaller than or equal to the one of $\BH^*_\alpha$. More precisely, since $\FDR(P,\BH^*_\alpha) = \alpha m_0/m$ under \eqref{indep} (see \cite{BY2001}), we have under \eqref{indep}, 
\begin{equation}
\FDR(P,\empBH_\alpha)\leq  \alpha m_0/m= \FDR(P,\BH^*_\alpha).\label{equ:conservative}
\end{equation}
In addition, the FDR control of $\empBH_\alpha$ holds under the more general condition \eqref{exch}. 
This is not the case for $\BH^*_\alpha$ that can violate the FDR control under that condition.
{Hence, Theorem~\ref{thFDR} puts forward }an additional robustness of $\empBH_\alpha$ w.r.t. the negative dependence, which is not enjoyed by $\BH^*_\alpha$. We provide an example below, see Figure~\ref{fig:side} for an illustration. 


 \begin{figure}[h!]
 \vspace{-1mm}
\begin{center}
\includegraphics[scale=0.6]{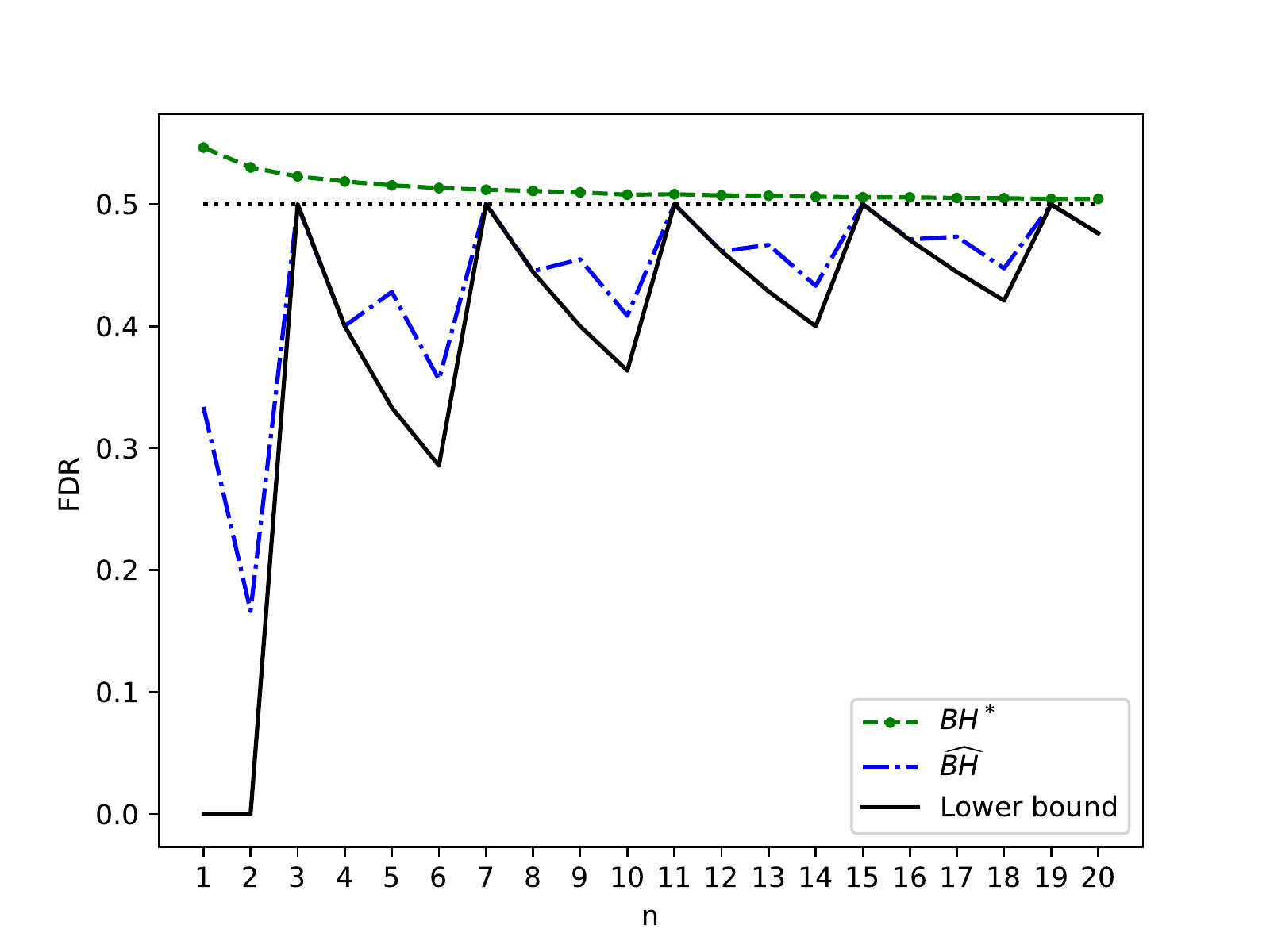}
\end{center}
\vspace{-5mm}
 \caption{\label{fig:side} {FDR of $\empBH$ and $\BH^*$ as a function of $n$, the size of the null training sample $Y$,  in a maximal negative dependent case (see Example~\ref{sec:robutnegative} below for more details). 
 The nominal level is $\alpha=0.5$ (horizontal dashed line). Here, $n$ ranges from $0$ to $20$, $m_0=m=2$,  and the FDR are evaluated with $10^6$ simulations. The FDR lower bound delineated in Theorem~\ref{thFDR} is also displayed. }}
 \end{figure}

 \begin{example}[Gaussian with maximal negative correlation] \label{sec:robutnegative}
 Assume that $Z=(Y,X)$ is a centered Gaussian vector with equicorrelation $\rho<0$ and variances equal to $1$. Classically, since the length of $Z$ is $n+m$, the condition $\rho\geq -1/(n+m-1)$ is necessary to provide that the $(n+m)\times (n+m)$ $\rho$-equicorrelated  matrix (that is, with diagonal $1$ and off-diagonal element $\rho$) is non-negative. For instance, the maximal negatively correlated case $\rho=-1/(n+m-1)$ can be easily realized as 
 $Z=(1+1/(n+m-1))^{1/2} (W_i-\ol{W})_{1\leq i\leq n+m},$ with $W_i$, $1\leq i\leq n+m$, i.i.d. $\mathcal{N}(0,1)$ and $\ol{W}$ denoting the sample mean of the $W_i$'s, $1\leq i\leq n+m$. For this specific distribution $P$ of $Z$, we have $P_0=P_0(P)=\mathcal{N}(0,1)$  and $\cH_0=\{1,\dots,m\}$. 
Also, Assumption \eqref{exch} is satisfied so that $\empBH_\alpha$ controls the FDR at level $\alpha$ (with equality when $\alpha (n+1)/m$ is an integer). On the other hand, it is well known that $\BH^*_\alpha$ has an FDR above $\alpha$ in that case (see also Figure~\ref{fig:side}). Additional illustrations are given in Section~\ref{sec:negative} in the numerical experiments. 
This example is also the starting point of the randomized BH procedure developed in Appendix~\ref{sec:randBH}.
 \end{example}

\begin{remark}\label{rempreviousFDR}
Since the first version of this work, earlier occurrences of the upper-bound proved in Theorem~\ref{thFDR} have been reported to us (our work has been developed independently): first, it has been proved under assumption \eqref{exch} in the work of \cite{weinstein2017power}  by using the same martingale as ours (in a different context). Second, the upper-bound is a consequence of the work of \cite{bates2021testing} who showed that the $\hat{p}$-values, despite their intricate structure, are positively regressively dependent on each one of the subset (PRDS). This is proved under the stronger assumption \eqref{indep}. 
\end{remark}

 \begin{remark}
 When $\alpha (n+1)/m$ is an integer, we can easily check that  $\empBH_\alpha$ coincide with $\widetilde{\BH}_\alpha$, the BH procedure applied to the naive, unbiased, $\tilde{p}$-values defined by \eqref{equptilde}. Hence, Theorem~\ref{thFDR} implies that  $\FDR(P,\widetilde{\BH}_\alpha)= \alpha m_0/m$ in that case (under \eqref{exch}). This shows that, perhaps surprisingly, the naive way to build empirical $p$-values eventually leads to a correct FDR control for such values of $n$. Simulations will show that this is not the necessarily the case for other values of $n$, see Section~\ref{sec:num}. 
 \end{remark}

\section{Power result}\label{powerresult}
\label{sec:powerrresult}

Section~\ref{sec:FDR} showed that $\empBH_\alpha$ has an FDR smaller than or equal to the one of the oracle $\BH^*_\alpha$ under \eqref{indep}, see \eqref{equ:conservative}. Now, an important concern is to check whether the {\it power} of $\empBH_\alpha$ is comparable to the one of $\BH^*_\alpha$. 
In this section, we explore this issue under Assumption~\eqref{indep} and the power comparison is established by comparing the true discovery proportions \eqref{equTDRTDP} of $\empBH_\alpha$ and $\BH^*_{\alpha'}$, for $\alpha'$ slightly below $\alpha$.
In a nutshell, we establish that the TDP of $\empBH_\alpha$ is larger than the one of $\BH^*_{\alpha'}$ with a probability tending to $1$, for any model parameter, when $n/m$ is large (Section~\ref{powerempBH}), while we show that it is not true when $n/m$ is small (Section~\ref{sec:impos}). Together, this means that the boundary achieved by the procedure $\empBH_\alpha$ is $n\asymp m/\alpha$. We then present the case of particular, more favorable distributions, for which at least  $k$ alternatives are ``detectable" (Section~\ref{sec:morefavorable}). In that case, the boundary achieved by $\empBH_\alpha$ is shown to be $n k \asymp m/\alpha$.

To state our results, let us finally introduce an additional notation: let 
\begin{align}\label{equ:pnm}
\mathcal{P}_{n,m}=\left\{P=P_0^{\otimes n}\otimes \bigotimes_{i=1}^m P_i\::\:  \mbox{$P_i$  continuous distribution on $\R$, $0\leq i\leq m$} \right\}.
\end{align}
Since  \eqref{indep} is true, $P$ belongs to this class $\mathcal{P}_{n,m}$ in the semi-supervised setting presented in Section~\ref{sec:setting}, which can thus be considered as the parameter set of the model under that assumption.
In addition, since we look at power results, we are going to focus on distributions in $\mathcal{P}_{n,m}$ with at least one true alternative. We denote 
\begin{align}\label{equ:anm}
\mathcal{A}_{n,m}=\{P\in\mathcal{P}_{n,m}\::\: m_1(P)\geq 1\}
\end{align}
 the corresponding set.

\subsection{Upper bound}\label{powerempBH}

The following result shows that, under \eqref{indep}, when $n\geq \gamma m$ with $\gamma$ large enough,  the semi-supervised BH procedure at level $\alpha$ rejects at least all null hypotheses rejected by the oracle BH procedure at level $\alpha'=\alpha(1-\eta)$, with high probability and with $\eta$ small.

\begin{proposition}\label{prop:power}
Recall   $\mathcal{A}_{n,m}$ \eqref{equ:anm}.
Let $\alpha, \gamma ,\eta \in (0,1)$ and let 
\begin{equation}\label{equgammaub}
\gamma^*(\alpha,\eta)=\alpha^{-1}\eta^{-2}(1+\eta)28\log(2)>0.
\end{equation}
 Then,  for all $n,m\geq 1$ with $n\geq \gamma m$, for all $P\in \mathcal{A}_{n,m}$, we have
\begin{equation}\label{equ:power}
\P_{Z\sim P}(\BH^*_{\alpha(1-\eta)}\subset \empBH_{\alpha})\geq 
1 -  { (1/2)^{3\gamma/\gamma^*(\alpha,\eta)-1}}.
\end{equation}
 In particular, for all $n,m\geq 1$ with $n\geq \gamma m$,
$$
\sup_{P\in \mathcal{A}_{n,m}}\{\P_{Z\sim P}(\TDP(P,\BH^*_{\alpha(1-\eta)})> \TDP(P,\empBH_{\alpha }))\}\leq 
{ (1/2)^{3\gamma/\gamma^*(\alpha,\eta)-1}}
$$
\end{proposition}

Proposition~\ref{prop:power} is proved in Section~\ref{proof:prop:power}.
It is based on a concentration argument of the empirical c.d.f. of the $Y_i$'s, which  relies on the independence assumption between the $Y_i$'s. 
Note that the bound \eqref{equ:power} is only informative if 
 $\gamma> \gamma^*(\alpha,\eta)/3$; otherwise the right-hand side of \eqref{equ:power} is non-positive and the bound is silent. When $\gamma
\geq  \gamma^*(\alpha,\eta)$, the probability is larger than or equal to $3/4$.
Taking  $\gamma$ much larger than $\gamma^*(\alpha,\eta)$ makes the probability arbitrarily close to $1$.

\subsection{Lower bound}\label{sec:impos}

The previous section shows that the power of $\empBH$ is close to the one of the oracle BH procedure provided that $n/m$ is sufficiently large. 
We can legitimately ask whether this condition is necessary.
The following result addresses this point.

\begin{proposition}\label{thm:lb0}
Recall  $\mathcal{A}_{n,m}$ \eqref{equ:anm}.
Let $\alpha\in (0,1)$ and $\eta\in (0,1)$. Consider $n,m\geq 1$ with $n/m \leq  1/(4\alpha)$. Then  \begin{equation}\label{equ:thmlb0}
 \sup_{P\in \mathcal{A}_{n,m}}\{\P_{Z\sim P}(\TDP(P, \BH^*_{\alpha(1-\eta)})>\TDP(P,\empBH_{\alpha}) )\}\geq 
 1-2\alpha.
 \end{equation}
\end{proposition}

Proposition~\ref{thm:lb0}  is proved in Section~\ref{proofthm:lb0}. It is a consequence of the fact that all $\hat{p}$-values are larger than $1/(n+1)$ (see \eqref{equemppvalue}), while $\empBH$ controls the FDR (Theorem~\ref{thFDR}).

Putting together, Propositions~\ref{prop:power}~and~\ref{thm:lb0} establish that the semi-supervised BH procedure achieves the boundary  $n\asymp m/\alpha$: for $n \leq  m/(4\alpha)$, there exists a configuration $P\in \mathcal{A}_{n,m}$ such that 
the power of $\empBH_{\alpha}$ is less than the one of the oracle $\BH^*_{\alpha(1-\eta)}$ (with probability at least $1-2\alpha$), while for $n\gg m/\alpha$ all configurations $P\in \mathcal{A}_{n,m}$ are such that 
the power of $\empBH_{\alpha}$ is larger than the one of the oracle  (with probability arbitrarily close to $1$).

\subsection{Refinement to more favorable distributions}\label{sec:morefavorable}

If there are enough alternatives, with enough signal strength, we show here that the  boundary achieved by $\empBH$ can be  much better than $m\asymp n/\alpha$. 
We extend for this Proposition~\ref{prop:power} and Proposition~\ref{thm:lb0} to a specific set of ``more favorable'' distributions.

For $\alpha\in (0,1)$, $n,m\geq 1$ and $1\leq k\leq m$, consider the subset of $\mathcal{A}_{n,m}$ given by
\begin{align*}
&\mathcal{A}_{n,m,k,\alpha,\beta}= \left\{P\in \mathcal{A}_{n,m}\::\: m_1(P)\geq k,\:\P_{Z\sim P}\left( 
| \cH_1(P)\cap \BH^*_{\alpha/2}|
\leq k-1\right)\leq \beta \right\}.
\end{align*}
In words, $\mathcal{A}_{n,m,k,\alpha,\beta}$ is the set of distributions such that at least $k$ null hypotheses are false while the probability that the procedure $\BH^*_{\alpha/2}$ makes at most $k-1$ number of true discoveries  is smaller than $\beta$. From an intuitive point of view, this means that the distribution contains at least $k$ ``detectable'' alternatives, in the sense that they are detectable with large probability by the oracle itself (at level $\alpha/2$). 

Now, the idea is that for a distribution $P\in \mathcal{A}_{n,m,k,\alpha,\beta}$, the threshold of the oracle procedure  is at least $\alpha k/m$ with large probability, so that the precision $1/n$ of the $\hat{p}$-values is enough to mimic the power of the oracle BH if and only if $1/n \ll \alpha k/m$, that is, $n k\gg m$. The following result proves that this informal argument is correct.

\begin{proposition}\label{extension}
Let $\alpha\in (0,1)$, $\eta\in (0,1/2)$ and $\beta\in (0,1)$. Then the following holds for $n,m\geq 1$ and $1\leq k\leq m$:
\begin{itemize}
\item[(i)] if $n k/m\geq \gamma $, for some $\gamma>0$,
$$
\sup_{P\in \mathcal{A}_{n,m,k,\alpha,\beta}}\{\P_{Z\sim P}(\TDP(P,\BH^*_{\alpha(1-\eta)})>  \TDP(P,\empBH_{\alpha }))\}\leq \beta+
{ (1/2)^{3\gamma/\gamma^*(\alpha,\eta)-1}}
$$
where $\gamma^*(\alpha,\eta)$ is given by \eqref{equgammaub}.
\item[(ii)] if $n k/m \leq  1/(4\alpha)$,
\begin{equation*}
 \sup_{P\in \mathcal{A}_{n,m,k,\alpha,\beta}}\{\P_{Z\sim P}(\TDP(P,\BH^*_{\alpha(1-\eta)})> \TDP(P,\empBH_{\alpha }))\}\geq 
 1-\beta-2\alpha.
 \end{equation*}
\end{itemize}
\end{proposition}

Proposition~\ref{extension} is proved in Section~\ref{sec:proofextension}.
Point (i) above is an upper-bound: in particular, it shows that having $n k/m\geq \gamma^*(\alpha,\eta) $ is enough for $\empBH_{\alpha }$ to mimic the power of the oracle $\BH^*_{\alpha(1-\eta)}$ with probability at least $1-\beta-1/4$ when the underlying distribution belongs to the set $ \mathcal{A}_{n,m,k,\alpha,\beta}$. 
Interestingly, the condition $n k/m\geq \gamma^*(\alpha,\eta) $ can be much weaker than the previous condition $n/m\geq \gamma^*(\alpha,\eta) $ when $k$ gets large. 

Point (ii) is a lower-bound showing that the order given in the upper-bound is correct. Together, (i) and (ii) ensure that the boundary achieved by $\empBH_{\alpha }$ is $nk\asymp m$ on the distribution set $\mathcal{A}_{n,m,k,\alpha,\beta}$. In addition, when $\alpha$ gets small and $1/\alpha$ cannot be considered as a constant, our result is able to track the dependence in $\alpha$; since $\gamma^*(\alpha,\eta)$ is of order $1/\alpha$ (see \eqref{equgammaub}), the boundary reads $nk\asymp m/\alpha$. This boundary turns out to be an accurate ``rule of thumb'' in the numerical experiments of Section~\ref{sec:num}.

\section{Optimality}\label{sec:optimality}

For a fixed level $\alpha$, the previous results show that the semi-supervised BH procedure $\empBH_{\alpha }$  mimics the oracle BH procedure when $n\gg m$ both in terms of FDR (Theorem~\ref{thFDR}) and power (Proposition~\ref{prop:power}).
However, when $n\ll m$, while $\empBH_{\alpha }$ still controls the FDR, it looses the power property (Proposition~\ref{thm:lb0}).  Hence, it does not mimic the oracle in that regime. 
However, this does not exclude that a different procedure, that would use the data $Z$ more cleverly, might be able to mimic the oracle when $n\ll m$. 
In this section, we show that {\it no procedure} can mimic the oracle in that regime (Theorem~\ref{thm:lb}). This 
shows a general  phase transition to the problem of mimicking the oracle  (Corollary~\ref{cor-lb}) and establishes that $\empBH_{\alpha }$ achieves this transition, which thus delineates a kind of optimality enjoyed by the semi-supervised BH procedure.

\subsection{General lower bound}

Recall  $\mathcal{P}_{n,m}$ \eqref{equ:pnm} and  $\mathcal{A}_{n,m}$ \eqref{equ:anm}.
Taken together, Theorem~\ref{thFDR} and Proposition~\ref{prop:power} show that for any $\alpha,\eta,\gamma\in (0,1)$ the  procedure $R=\empBH_{\alpha }$ (as a sequence in $n,m\geq 1$) enjoys simultaneously the two following properties:
\begin{align}
\sup_{\substack{n,m\geq 1\\n\geq m\gamma}} \sup_{P\in \mathcal{P}_{n,m}}\{\FDR(P,R) - \FDR(P,\BH^*_\alpha)\}&\leq \delta_1;\label{mimic1}\\
\sup_{\substack{n,m\geq 1\\n\geq m\gamma}}\sup_{P\in \mathcal{A}_{n,m}}\P(\TDP(P,\BH^*_{\alpha(1-\eta)})> \TDP(P,R))&\leq \delta_2.\label{mimic2}
 \end{align}
for $\delta_1=0$ and $\delta_2=(1/2)^{3\gamma/\gamma^*(\alpha,\eta)-1}>0$. 
This quantifies how $\empBH_{\alpha }$ mimics the oracle $(\BH^*_\alpha)_{\alpha\in (0,1)}$ both in terms of FDR and power when $\gamma/\gamma^*(\alpha,\eta)$ grows. 

In the regime where $\gamma$ is too small, the following result shows that achieving simultaneously \eqref{mimic1}  and \eqref{mimic2} is not possible, and this for any procedure $R$ based only on the data $Z$.

\begin{theorem}\label{thm:lb}
Recall  $\mathcal{P}_{n,m}$ \eqref{equ:pnm} and  $\mathcal{A}_{n,m}$ \eqref{equ:anm}.
 Let $\alpha\in (0,1/4)$ and $\gamma,\eta\in (0,1)$ and let 
\begin{equation}\label{equgammalb}
\gamma_*(\alpha,\eta)=(1+(\alpha(1-\eta))^{-1/2})^{-3}/64,
\end{equation}
Consider $n,m\geq 1$ with $n\leq \gamma m$. Then for any procedure $R$ (based only on $Z$), either
\begin{equation}\label{equ:thmlb0}
\sup_{P\in \mathcal{P}_{n,m}}\{\FDR(P,R)- \FDR(\BH^*_{\alpha},R)\}\geq 1/2-\alpha
 \end{equation}
or 
 \begin{equation}\label{equ:thmlb}
 \sup_{P\in \mathcal{A}_{n,m}}\{\P_{Z\sim P}(\TDP(P,R)< \TDP(P, \BH^*_{\alpha(1-\eta)}))\}\geq 
 1/2- (1/4) (\gamma/\gamma_*(\alpha,\eta))^{1/3}.
 \end{equation}
As a result,
no procedure $R$ (as a sequence in $n,m\geq 1$) can satisfy simultaneously \eqref{mimic1}  and \eqref{mimic2} for $ \delta_1<1/2-\alpha$ and $\delta_2<1/2- (1/4) (\gamma/\gamma_*(\alpha,\eta))^{1/3}$.
\end{theorem}

The proof of Theorem~\ref{thm:lb} is given in Section~\ref{proof:thm:lb}. It relies on building two nearly indistinguishable configurations $Q_1,Q_2\in \mathcal{P}_{n,m}$ such that: either $\FDR(Q_1,R)$ is large, or with large probability under $Q_2$, $R$ makes no discovery while the oracle makes at least one correct discovery.
Note that the result in Theorem~\ref{thm:lb} is silent if $1/2- (1/4) (\gamma/\gamma_*(\alpha,\eta))^{1/3}\leq 0$, that is, $\gamma\geq 8\gamma_*(\alpha,\eta)$. Hence, Theorem~\ref{thm:lb} is only informative  whenever $\gamma< 8\gamma_*(\alpha,\eta)$. When $\gamma< \gamma_*(\alpha,\eta)$, the RHS of \eqref{equ:thmlb} is in addition strictly larger than $1/4$.

\subsection{Phase transition}\label{sec:boundary}

Let us elaborate further on the phase transition that we have put forward. 
To this end, we introduce the following definition.

\begin{definition}
For $\delta_1\in [0,1)$, $\delta_2\in [0,1)$, $\gamma>0$ and $\alpha,\eta \in (0,1)$, a procedure $R$ is said $(\delta_1,\delta_2)$-mimicking the oracle $(\BH^*_\alpha)_{\alpha\in (0,1)}$, for a training-to-test sample size at least $\gamma$, and a level $\alpha$ with relaxation $\eta$, in short $R$ is $MO(\gamma,\alpha,\eta,\delta_1,\delta_2)$, when \eqref{mimic1} and \eqref{mimic2} simultaneously hold for these values of $\delta_1,\delta_2,\alpha,\gamma,\eta$.
\end{definition}

According to this definition, Theorem~\ref{thFDR}, Proposition~\ref{prop:power} and Theorem~\ref{thm:lb} can be combined as follows:

\begin{corollary}\label{cor-lb}
Let $\alpha\in (0,1/4),$ $\eta\in (0,1)$, and consider 
$\gamma^*(\alpha,\eta)
$ defined by \eqref{equgammaub}
and
$\gamma_*(\alpha,\eta)$ defined by \eqref{equgammalb}. Then for any $\gamma>0$:
\begin{itemize}
\item[(i)] If $\gamma < \gamma_*(\alpha,\eta)$, then there exists no $MO(\gamma,\alpha,\eta,\delta_1,\delta_2)$ procedure  for any possible value of $ \delta_1,\delta_2\in (0,1/4]$. This is {even true for $\delta_1< 1/2-\alpha$ and $\delta_2<1/2-  (1/4) (\gamma/(\gamma_*(\alpha,\eta)))^{1/3})$};
\item[(ii)] If $\gamma\geq \gamma^*(\alpha,\eta)$ then there exists an $MO(\gamma,\alpha,\eta,\delta_1,\delta_2)$ procedure for some values of $\delta_1,\delta_2\in [0,1/4]$. This is {achieved by $\empBH_\alpha$, even with $\delta_1=0$ and $\delta_2=(1/2)^{3\gamma/\gamma^*(\alpha,\eta)-1}$.}
\end{itemize}
\end{corollary}

{This phase transition is illustrated in Figure~\ref{fig:phasetransition} in the introduction of the paper. The transition is provided under the slightly modified form $n=m/\alpha$ for comparison with Section~\ref{powerresult}. 
This also emphasizes that the impossibility result is a worst case analysis over the distribution $P\in \mathcal{P}_{n,m}$ (FDR) and $P\in \mathcal{A}_{n,m}$ (power) (suprema are taken in \eqref{equ:thmlb0} and \eqref{equ:thmlb}). In particular, under the more stringent assumption  $P\in \mathcal{A}_{n,m,k,\alpha,\beta}$, mimicking the oracle becomes already  possible whenever $n\gtrsim m/(\alpha k)$ (as reported in Figure~\ref{fig:phasetransition} for $k=3$ or $100$).
}

This general phase transition is in line with the recent results by \cite{roquain2020false}. 
Nevertheless, their setting is markedly different: it is unsupervised (no NTS) and the null distribution is assumed to belong to the Gaussian distribution family with unknown mean and variance. The phase transition found there was $\l \asymp m/\log(m)$ (with our notation) where $\l$ is a lower bound on the number of alternatives $m_1(P)$ (with no signal strength assumption). Here, the situation is notably different, with a boundary function of the length $n$ of the NTS.
The situation is also very different in terms of FDR control: the mimicking procedure $\empBH$ provides here an FDR control both above and below the transition, while such property is not possible in the setting of \cite{roquain2020false} (as proved in Corollary~3.3 therein).

\section{{Numerical illustrations}}\label{sec:num}

This section provides several numerical illustrations for the theoretical findings derived in Sections~\ref{sec:FDR}~and~\ref{powerresult}.

\subsection{Simulation setting} 

 While our experiments mostly focus on the two BH-type procedures $\empBH$ and $\BH^*$, we will also consider other competitors: $\widetilde{\BH}$, which is the BH procedure applied to the unbiased $\tilde{p}$-values defined by \eqref{equptilde} (Section~\ref{sec:pcp}) and the ``naive'' procedures $\empBY$ and $\empBHsplit$ described in Section~\ref{sec:relatedwork}.
 Also, for simplicity, the way to evaluate how the power of $\empBH$ mimics the one of $\BH^*$ slightly departs from our theoretical study: first, we compare $\empBH$ to the oracle $\BH^*$ taken at the same level $\alpha$ (say, $\eta=0$ with the notation of Section~\ref{powerresult}). This makes the power mimicking more challenging. Second, to stick with the standard way of comparing procedures (for $\BH^*$, $\empBH$ or { their } competitors), the considered power criterion is simply the TDR \eqref{equTDRTDP} (average of the TDP). 
Unless specified, the setting is Gaussian with a null distribution $P_0=\mathcal{N}(0,1)$ and an alternative $\mathcal{N}(\mu,1)$, for a given value of $\mu>0$. 
Across the sections below, we made various choice{s} of $n$, $m$ and of the sparsity $m_1$ (number of alternatives). We sometimes fix the level $\alpha$ to the (unusual large) value $0.5$ for better visibility of the curves and faster computation time, but the results scale accordingly for smaller values of $\alpha$. Finally, the FDR (resp. TDR) curves are here estimated by Monte-Carlo simulations. The plots show the estimates $\widehat{\FDR}$ and $\widehat{\TDR}$ with two error bars: one estimating the standard deviation of $\widehat{\FDR}$ (resp. $\widehat{\TDR}$) and one estimating the standard variation of $\FDP$ (resp. $\TDP$) (these two deviations being proportional).

\begin{figure}[h!]
    \centering
    \subfigure[]{\includegraphics[width=.51\textwidth]{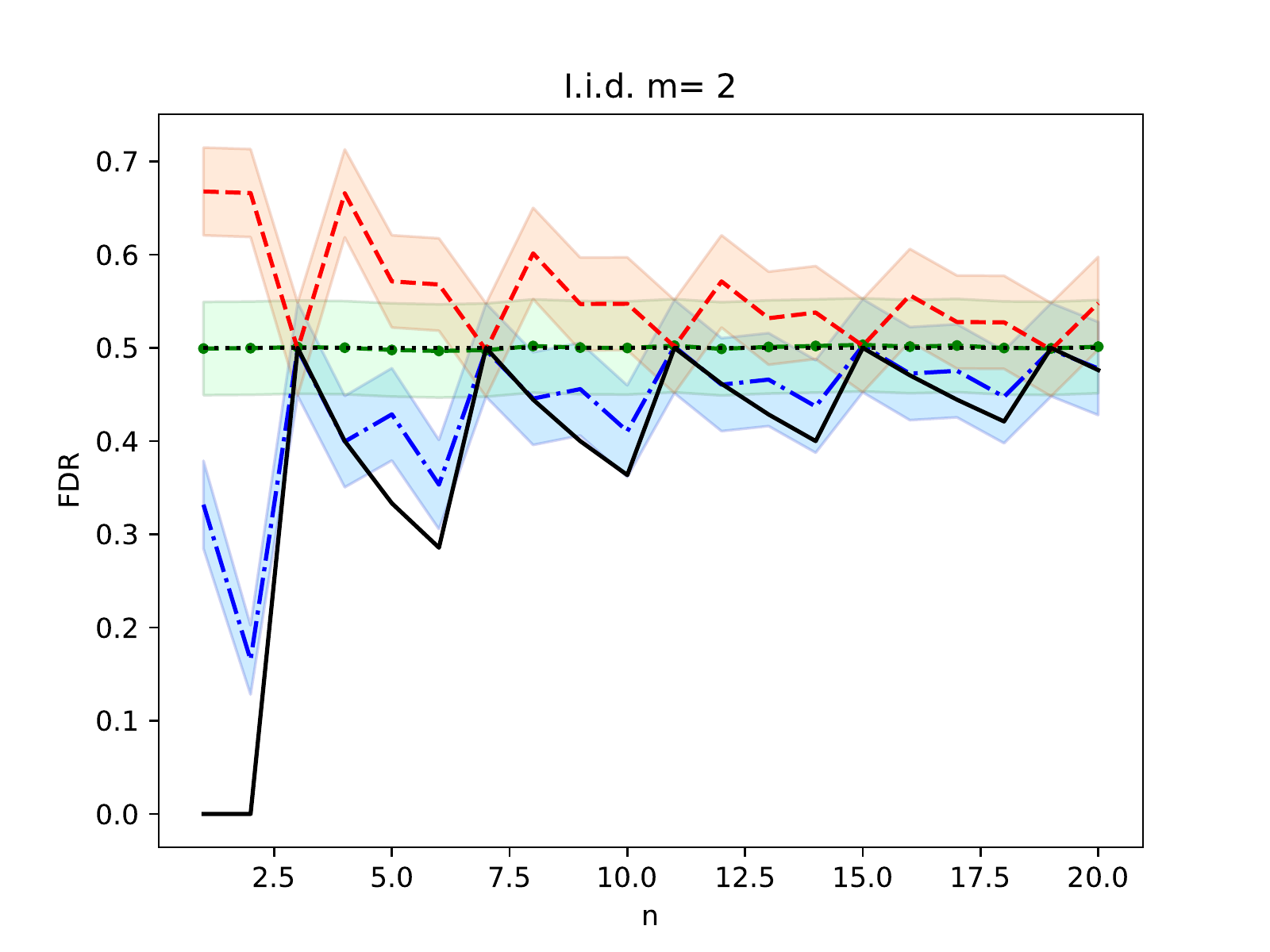}} \hspace{-9mm} \vspace{-3mm}
    \subfigure[]{\includegraphics[width=0.51\textwidth]{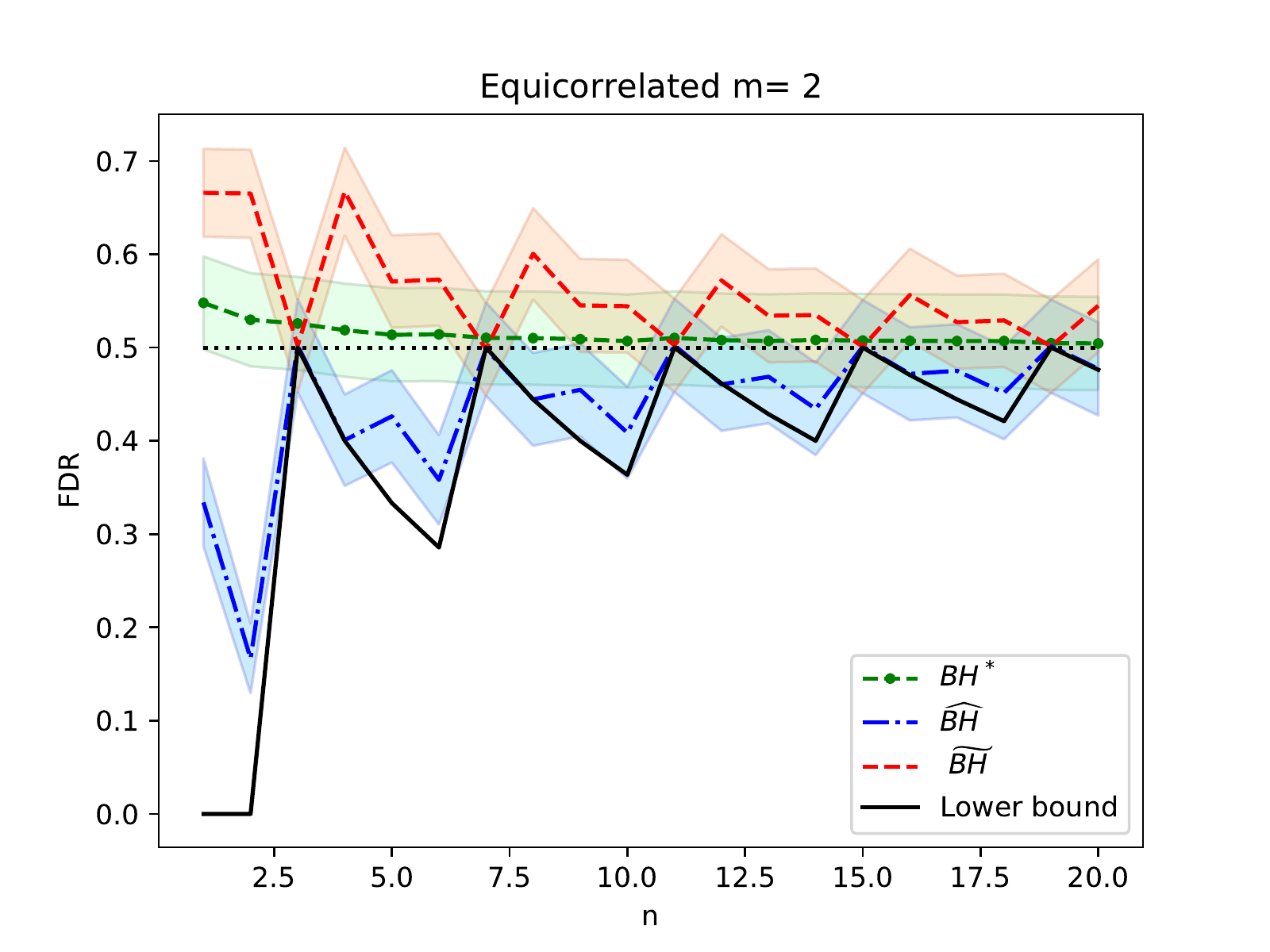}} \vspace{-3mm}
    \subfigure[]{\includegraphics[width=.51\textwidth]{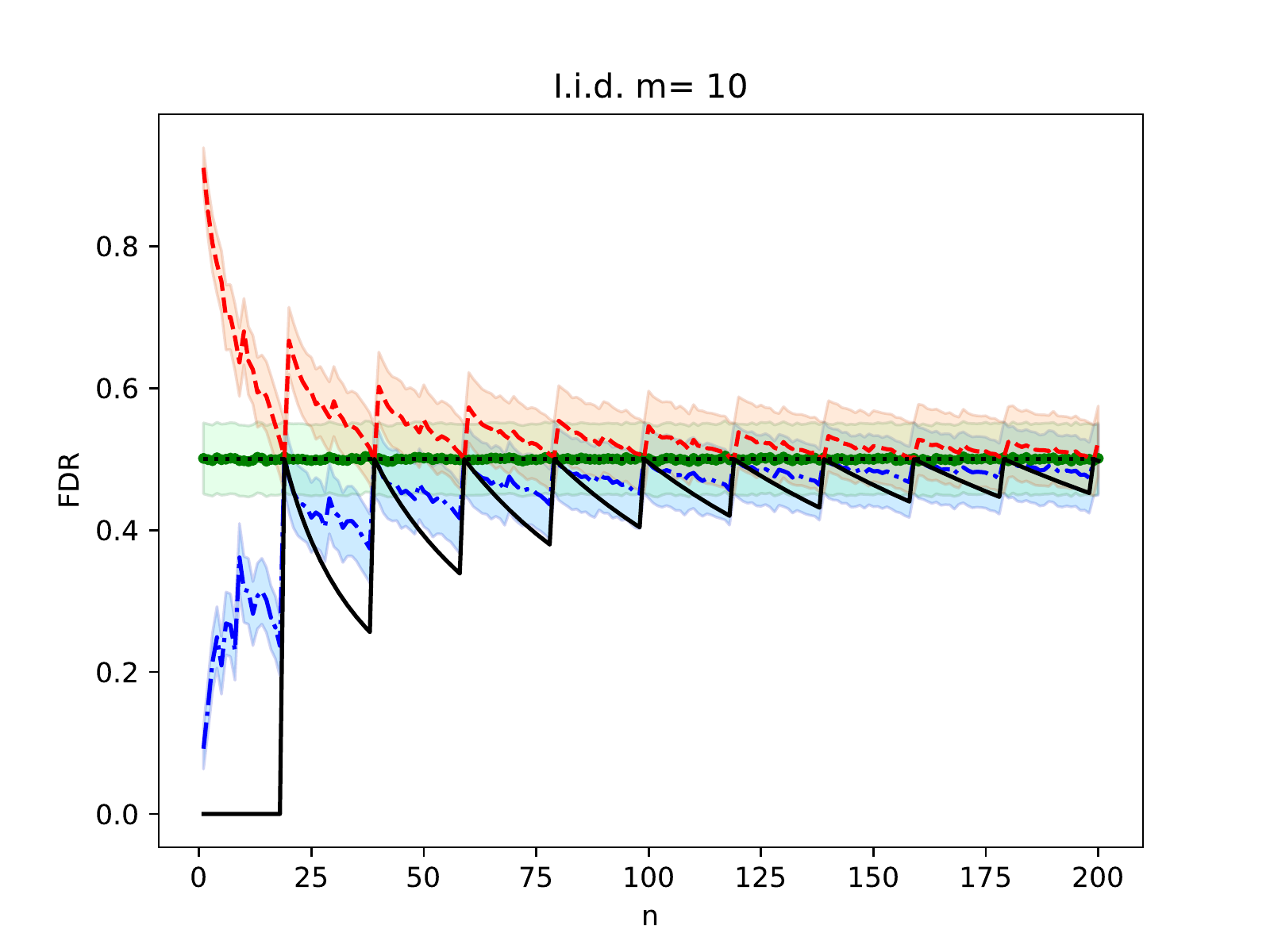}} \hspace{-9mm}
    \subfigure[]{\includegraphics[width=0.51\textwidth]{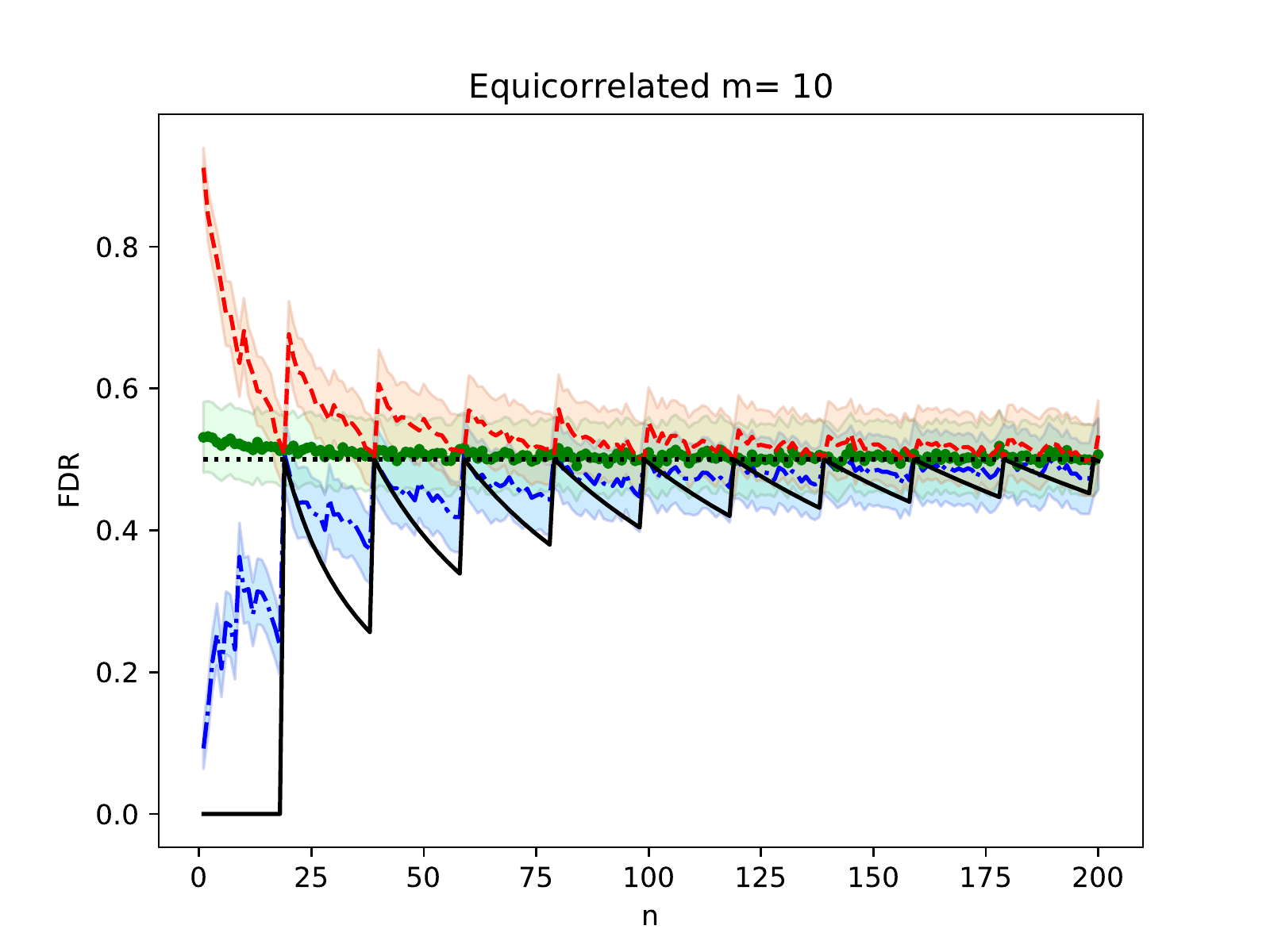}} 
     \caption{\label{fig:fdr}  FDR result in the case of { i.i.d. samples  (left column) and Gaussian negative equicorrelation  (right)}. The cases $m=2$ (top row), 
    $m=10$ (bottom row)  have been investigated with respectively $10^5$ and $10^4$ Monte Carlo simulations.  {The $2\sigma$ confidence intervals on the estimated FDR are not visible}. The standard deviation (divided by a factor $10$) of the FDP  is shown by shaded areas. The figure shows the results for $\BH^*$ (green), $\empBH$ (blue), $\widetilde{\BH}$ (red, see text), and the lower bound  of Theorem~\ref{thFDR}  (black). 
    } 
        \end{figure}

\subsection{FDR control under the full null} 
\label{sec:negative}

The first experiment concerns the case where $m_0=m$, which corresponds to the so-called ``full null'' configuration where there is no alternative. 
We consider two dependence framework: the independent case (all $Z_i$'s independent) and the negatively equicorrelated case described in Example~\ref{sec:robutnegative}. Recall that $\empBH$ is proved to control the FDR at level $\alpha$ in both cases  (Theorem~\ref{thFDR}), while $\BH^*$ is only proved to control the FDR at level $\alpha m_0/m=\alpha$ in the independent case. 
Also, $\widetilde{\BH}$ (BH procedure applied to the unbiased $\tilde{p}$-values defined in Section~\ref{sec:pcp}) is not proved to control the FDR since the $\tilde{p}$-values do not satisfy the super-uniformity property \eqref{equ:superunif}. 

Figure~\ref{fig:fdr} displays the obtained FDR curves for $\BH^*$ (green), $\empBH$ (blue) and $\widetilde{\BH}$ (red). The obtained results are consistent with our theoretical findings: negative correlations induce an FDR of $\BH^*$ slightly above the targeted level, although this effect tends to reduce when $n$ gets larger. This is because the negative correlation $\rho=-(m+n-1)^{-1}$ decreases (in absolute value) when $n$ grows.
As expected, $\empBH$ maintains the FDR control is any case. 
Meanwhile, $\widetilde{\BH}$ fails to control the FDR in any case, except for some values of $n$ where it has the same FDR value as $\empBH$. Hence, we shall discard $\widetilde{\BH}$ from our plots in the sequel. Interestingly, we also displayed the lower bound of Theorem~\ref{thFDR} in Figure~\ref{fig:fdr}: while it correctly lower bounds  the estimated FDR of $\empBH$ for any $n$, it illustrates that the FDR is exactly $\alpha$ for $n\in \{3,7,11,15,\dots\}$ as the theory establishes (the curves might also suggest that the lower bound is sharp for $n\in \{4,8,12,16,\dots\}$, which is not covered by our theory). Finally, note that these results are in expectation: as shown by the shaded areas, there can be  large variations for particular samples. This {is} inherent to the BH procedure when the number of discoveries is not large.

 \begin{figure}[h!]
  \vspace{-1mm}
    \centering
    {\includegraphics[width=.5\textwidth]{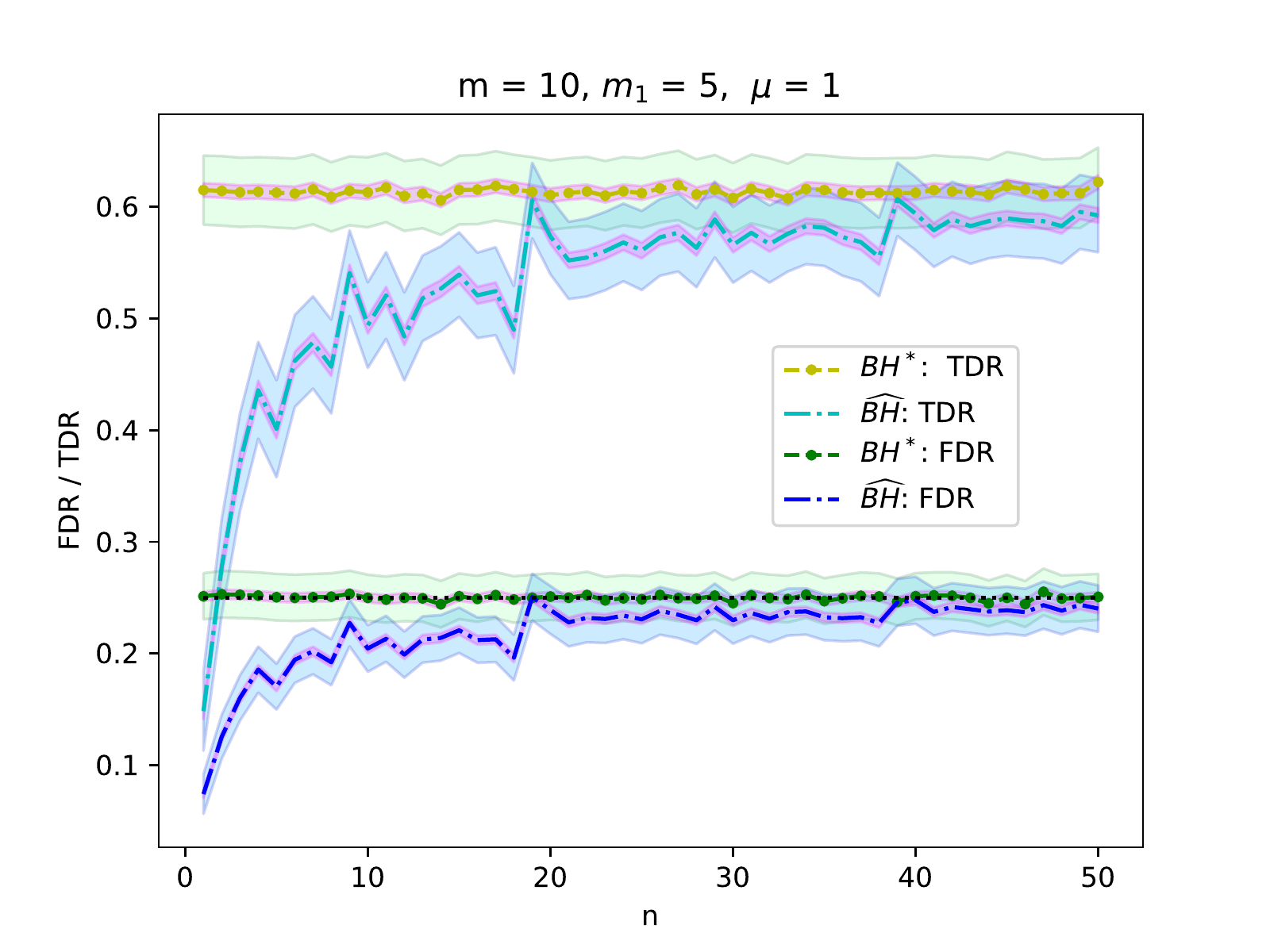}} \hspace{-8mm} 
  \includegraphics[width=.5\textwidth]{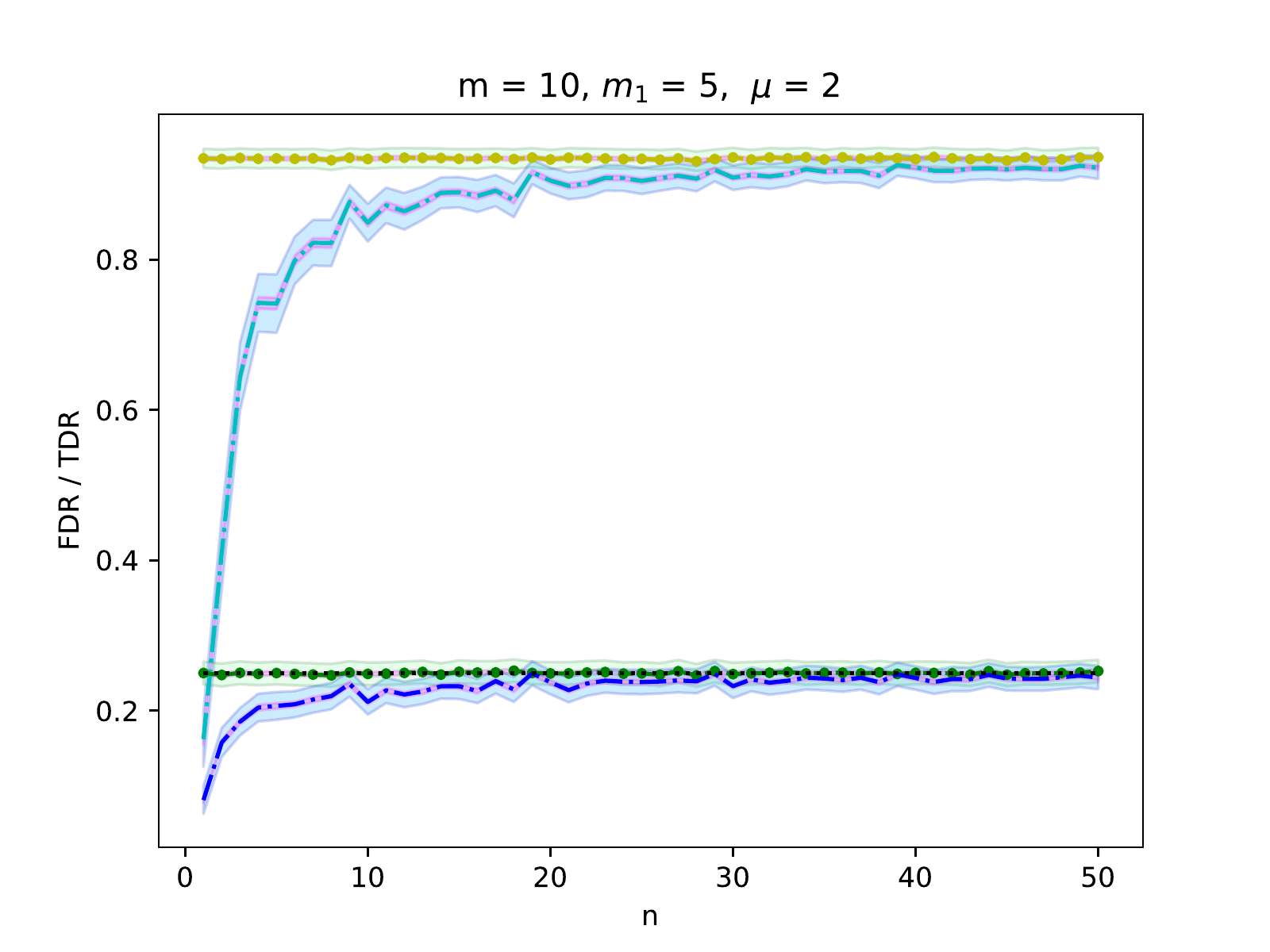} \\
      {\includegraphics[width=.5\textwidth]{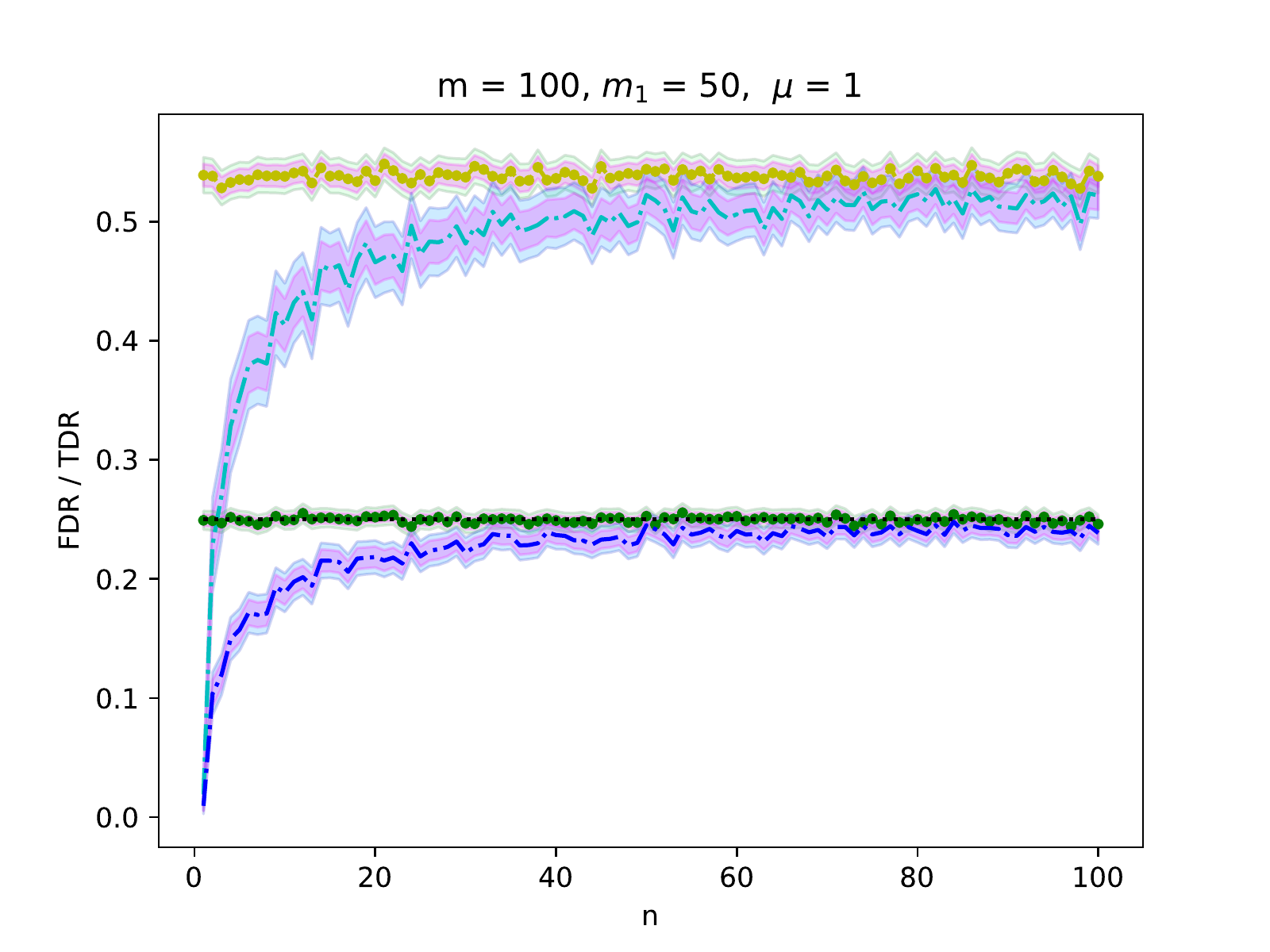}}\hspace{-8mm} 
 \includegraphics[width=.5\textwidth]{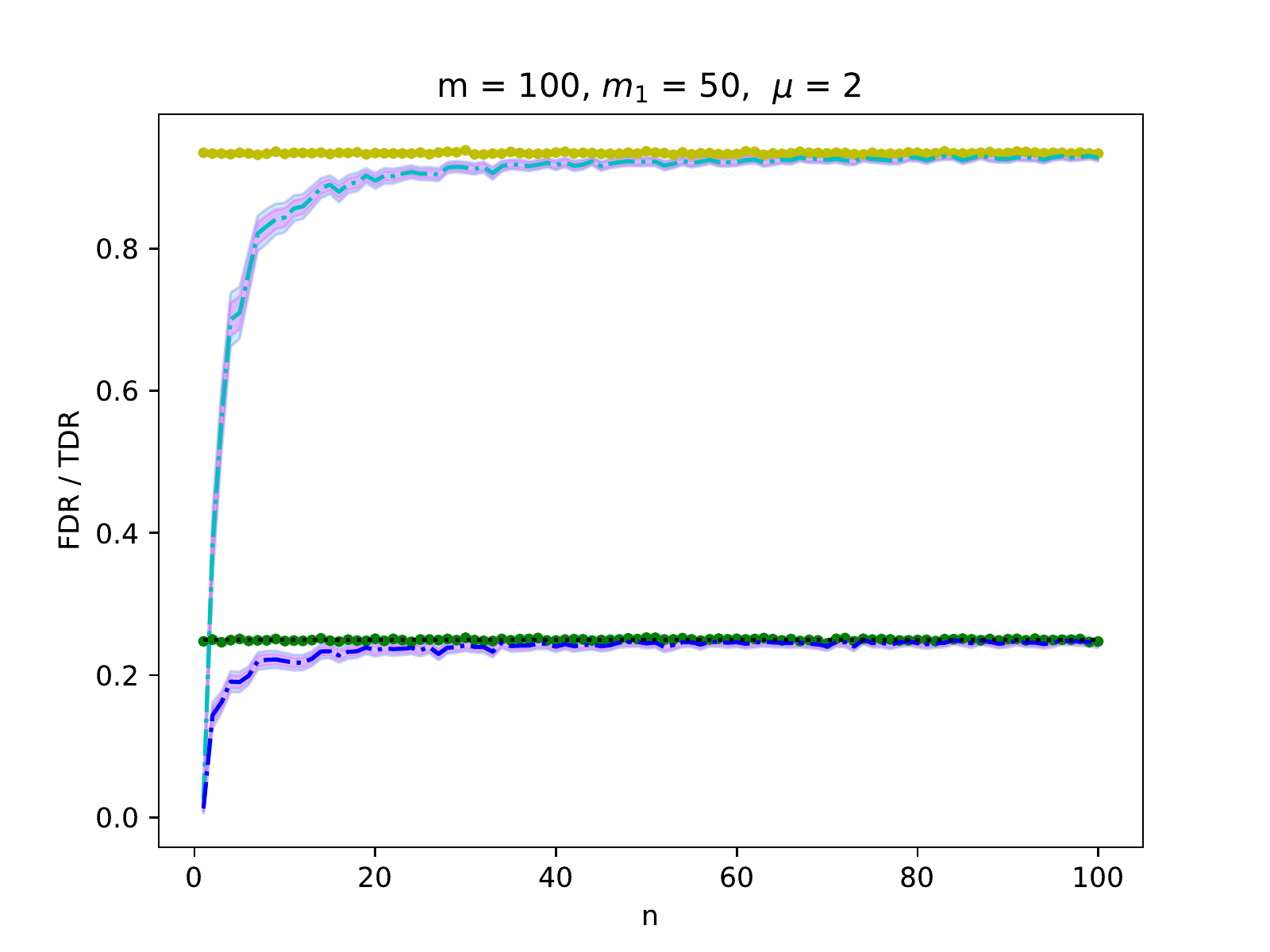} \vspace{-3mm}
     \caption{\label{fig:puis1} 
    FDR and TDR results for the dense case: ${m_1}=\frac{m}{2}$, with $\mu =1$ (left column) and $\mu =2$ (right column). The number of tests $m$ equals  $10$ in the top row
    and $100$ in the bottom row. The number of Monte Carlo simulations used for estimating the FDR and TDR is $10^4$ (top row) and $10^3$ (bottom  row).  The $2\sigma$ confidence interval on the estimated FDR and TDR is plotted in magenta. In all plots the standard deviation (divided by $10$) of the FDP and TDP are  shown in shaded green for $\BH^*$ and shaded blue for $\empBH$. }
    
      \end{figure}

\begin{figure}[h!]
  \vspace{-1mm}
    \centering
  {\includegraphics[width=.51\textwidth]{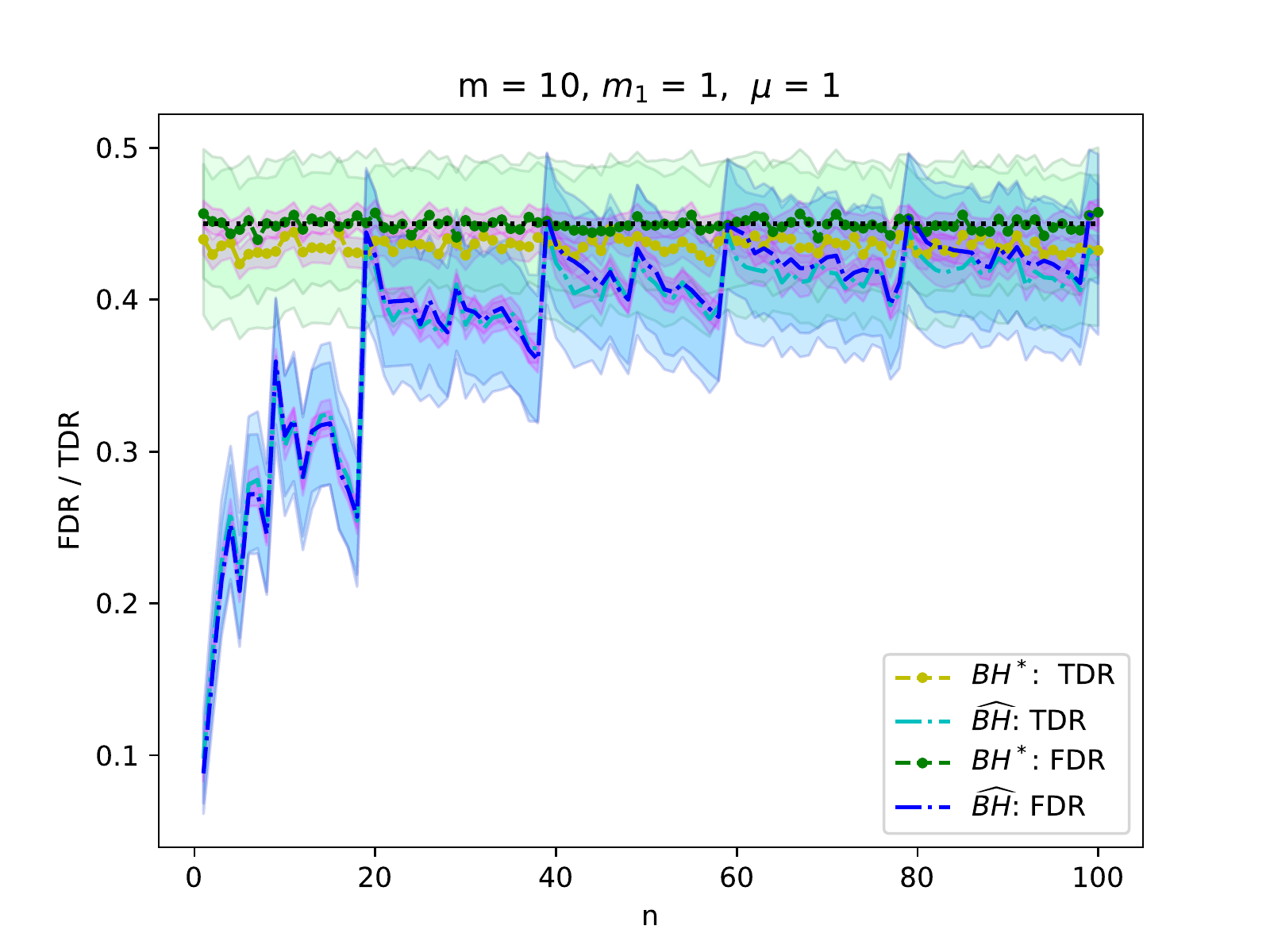}}\hspace{-8mm} 
 \includegraphics[width=.51\textwidth]{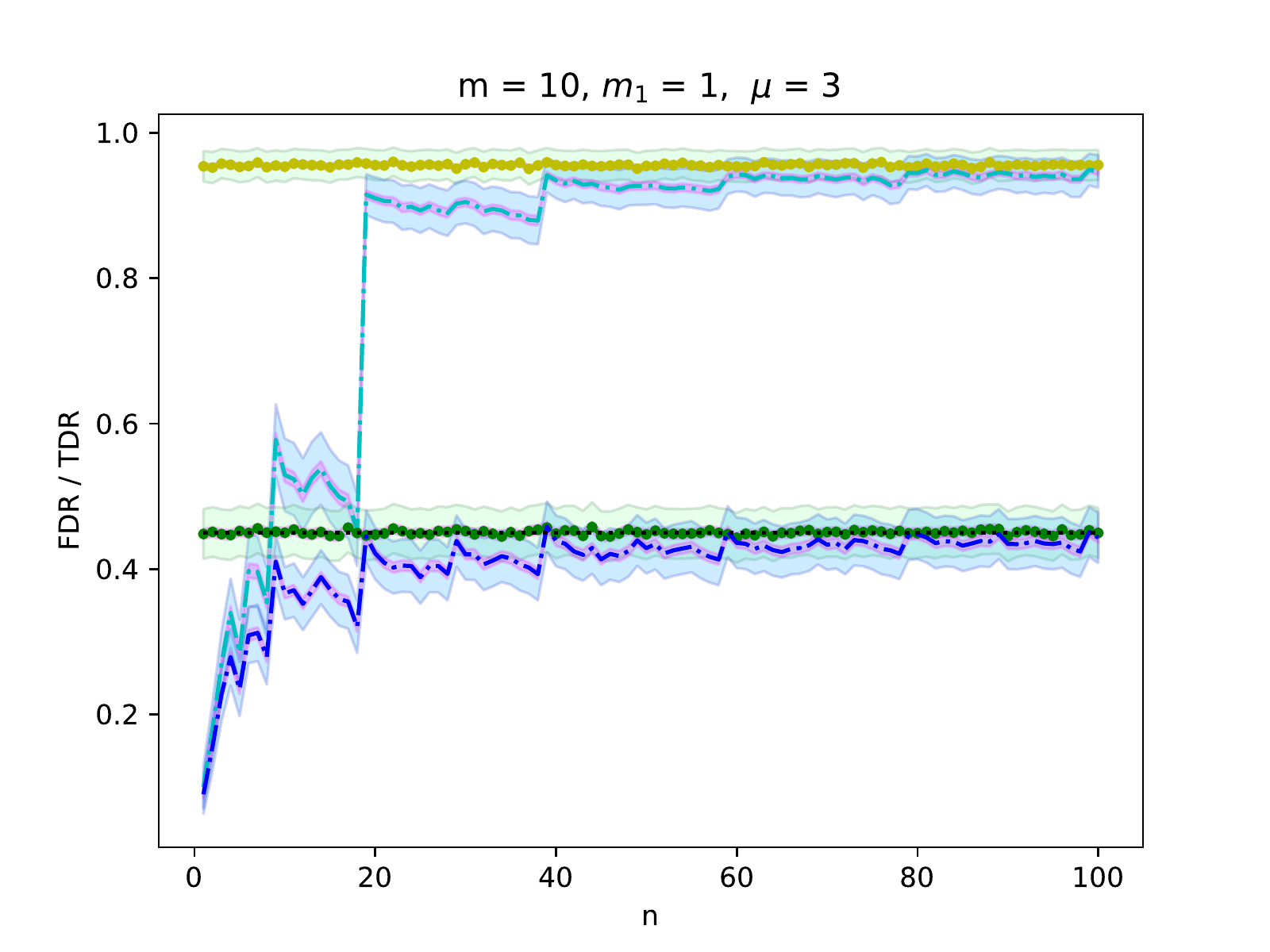} \\
          {\includegraphics[width=.51\textwidth]{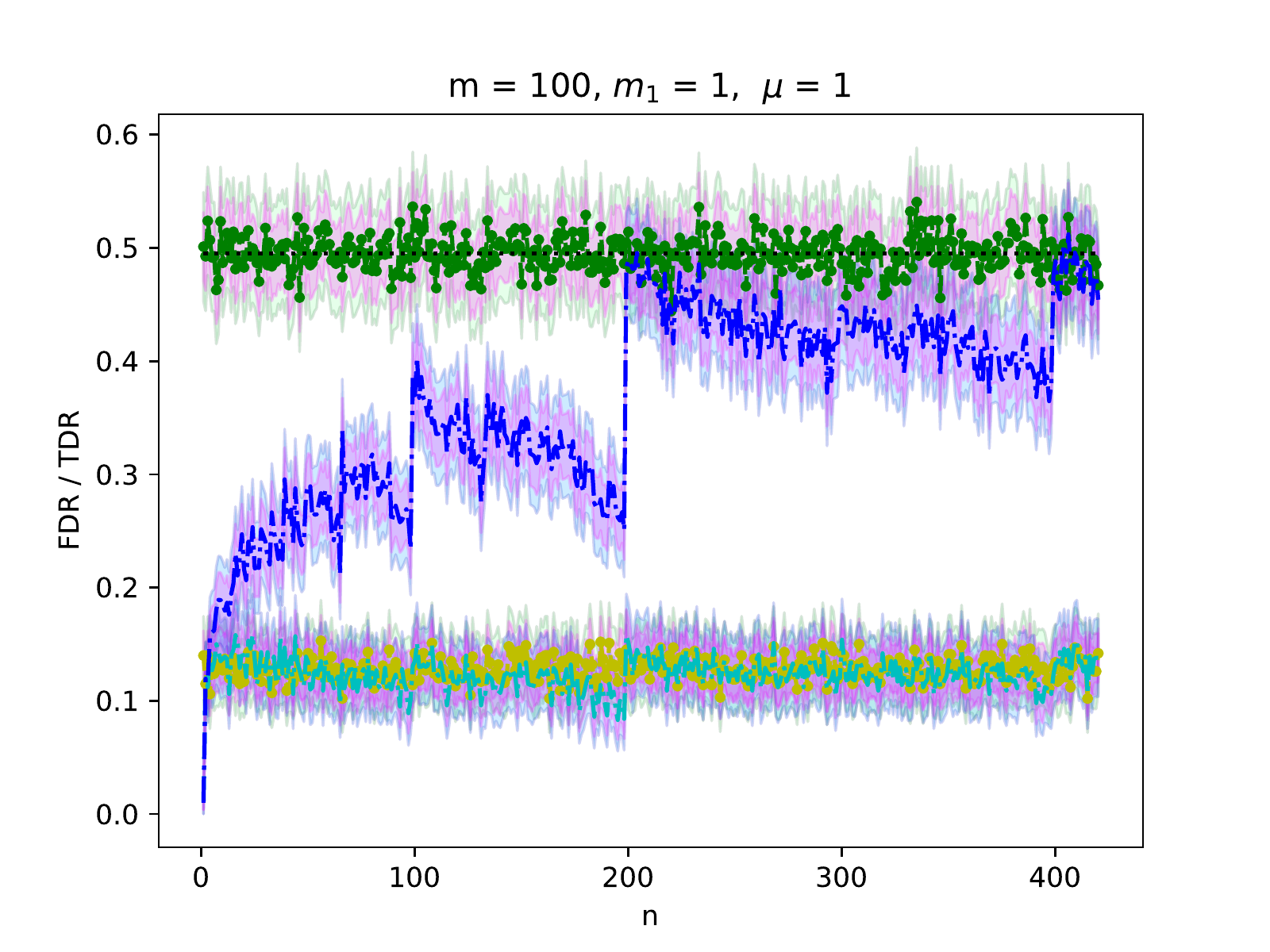} \hspace{-8mm} 
           \includegraphics[width=.51\textwidth]{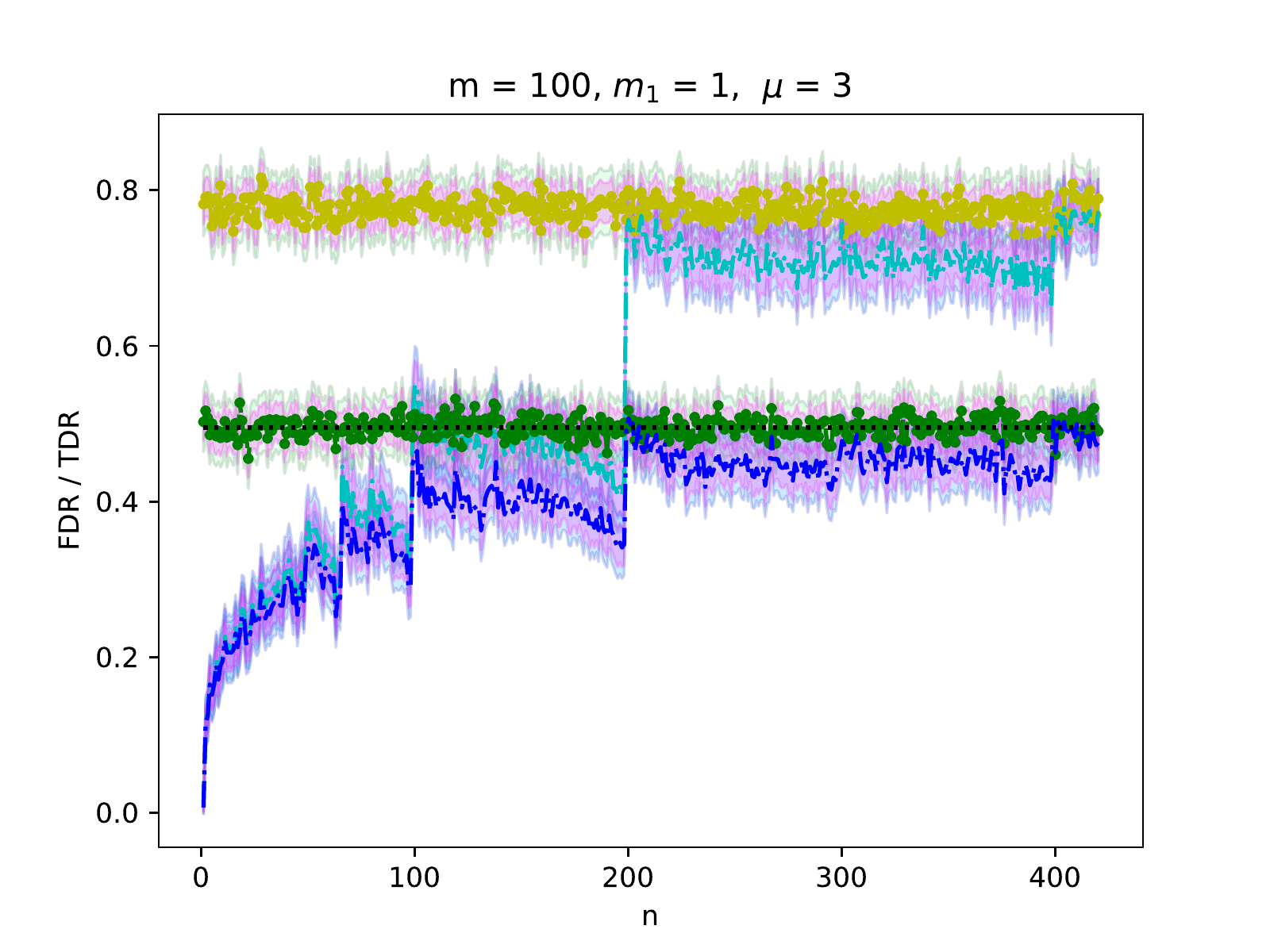}} \\
    \caption{\label{fig:puis2} FDR and TDR results for the sparse case: ${m_1}=1$, with $\mu =1$ (left column) and $\mu =3$ (right column). The number of tests $m$ equals  $10$ for the top row
    and $100$ for the bottom row. The number of Monte Carlo simulations used for estimating the FDR and TDR is $10^4$ (top row) and $10^3$ (bottom  row).  The $2\sigma$ confidence interval on the estimated FDR and TDR is plotted in magenta. In all plots the standard deviation (divided by $10$) of the FDP and TDP are  shown in shaded green for $\BH^*$ and shaded blue for $\empBH$.  }    
\end{figure}

\begin{figure}
  \vspace{-1mm}
    \centering
  {\includegraphics[width=.51\textwidth]{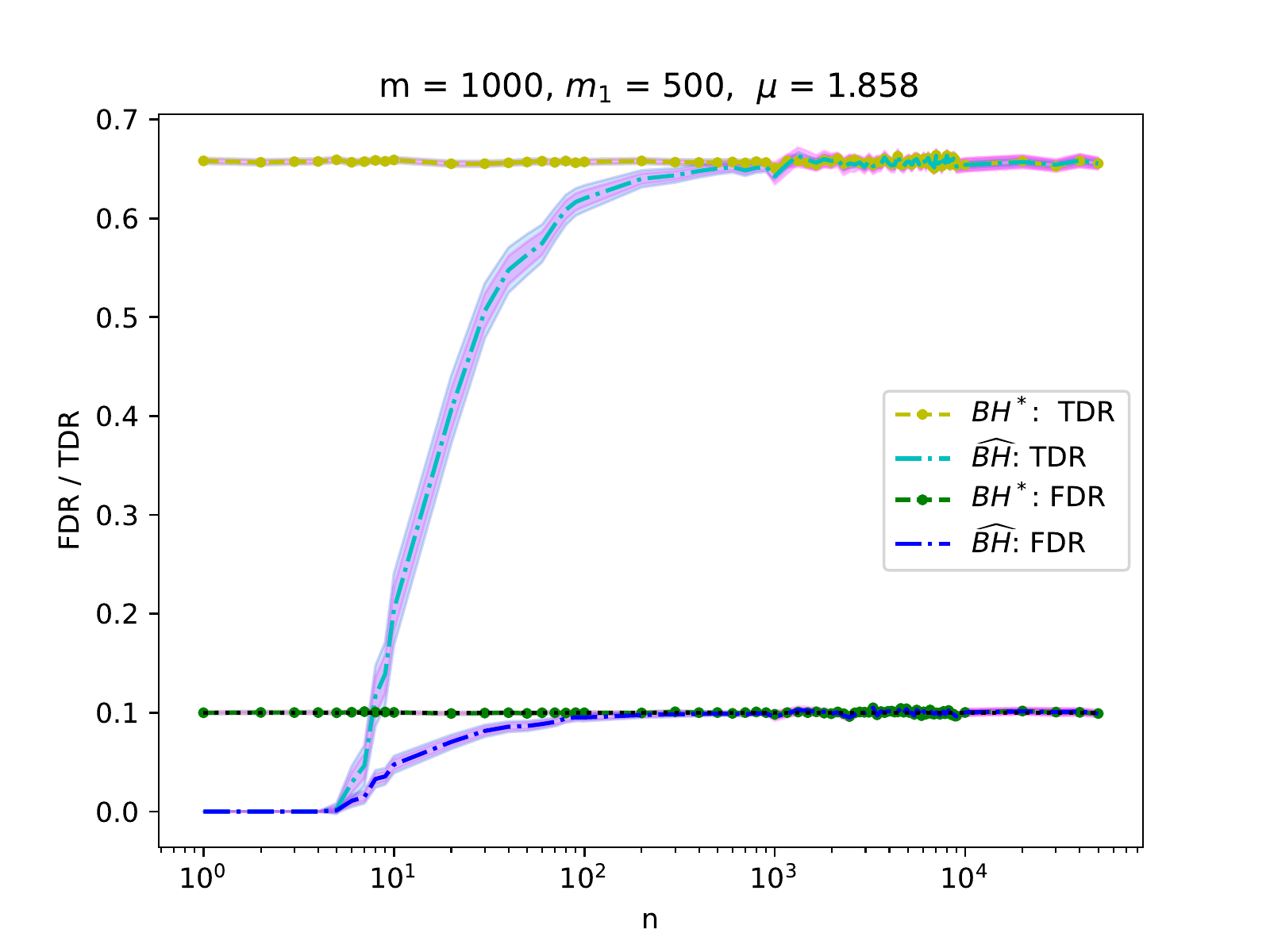}} \hspace{-8mm} 
 \includegraphics[width=.51\textwidth]{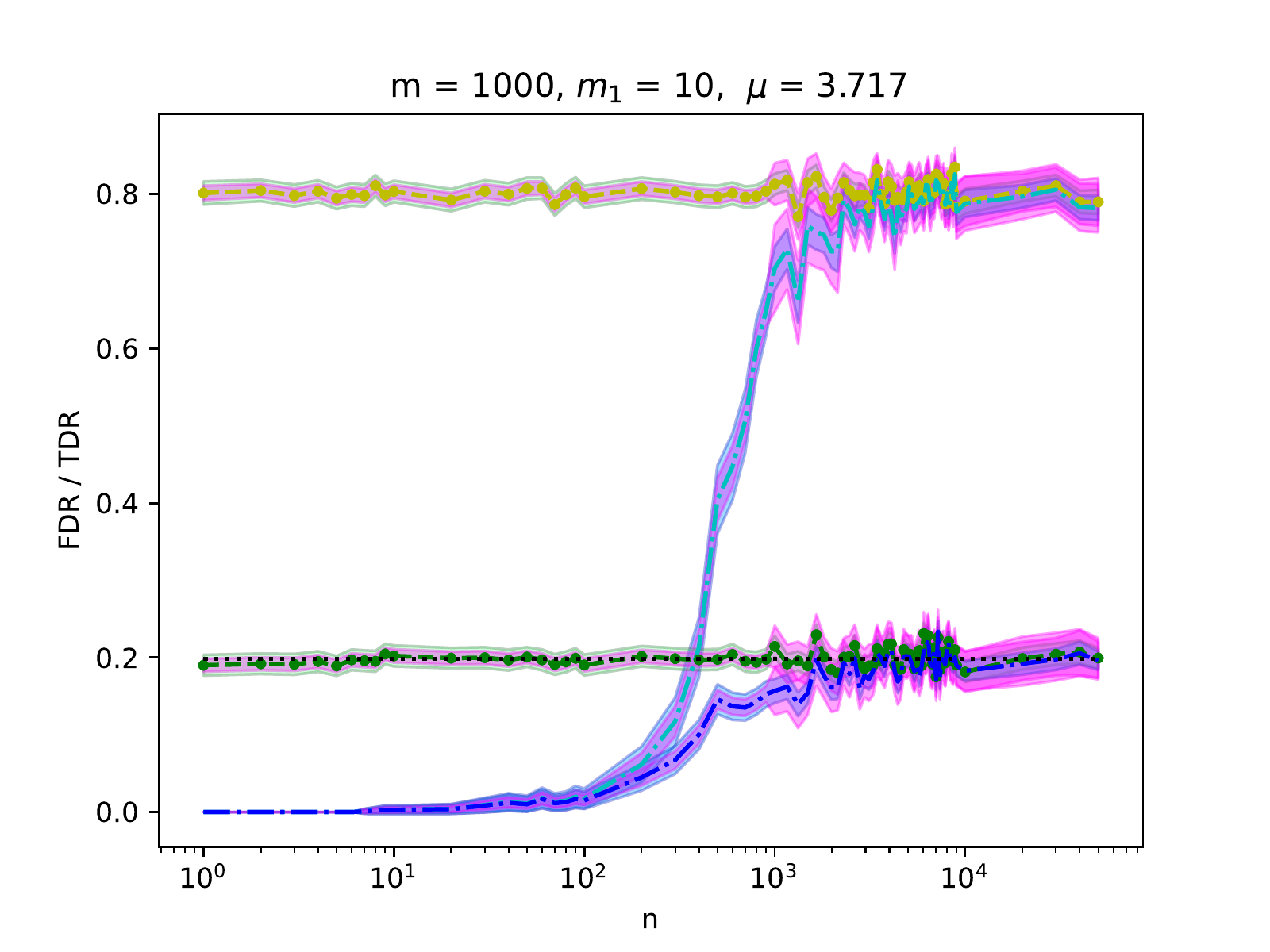} \\
    \caption{\label{fig:puis3} FDR and TDR results for $m=10^3$ as a function of $n$, for $\alpha = 0.2$. Left : dense case: $m_1=500$ and  $\mu =\frac{1}{2}\sqrt{2\log m}\approx 1.86$. Right: $m_1=10$ and  $\mu =\sqrt{2\log m}\approx 3.72$ (right). The number of Monte Carlo simulations used for estimating the FDR and TDR is $10^3$ for $n<10^3$ and   $10^2$ otherwise.  The $2\sigma$ confidence interval on the estimated FDR and TDR is plotted in magenta. In all plots the standard deviation (divided by $10$) of the FDP and TDP are  shown in shaded green for $\BH^*$ and shaded blue for $\empBH$.  }    
\end{figure}

\subsection{Power study }
\label{sec:puiis}

Figure~\ref{fig:puis1} compares the performances of the procedures $\BH^*$ (dark green and khaki) and $\empBH$ (dark blue and cyan) in terms of FDR (dark colors) and TDP (light colors) in the {dense} case where $m_1=m_0=\frac{m}{2}$, with $\mu =1$ (left column) and $\mu =2$ (right column).
Regarding the FDR first, the plots show that the FDR of $\empBH$
 tends to the oracle FDR (which is $0.25=\alpha \frac{m_0}{m} =\frac{\alpha}{2}$ here). {For a fixed value of $n$, the convergence is faster for smaller values of $m$}. {This is coherent with Theorem~\ref{thFDR}, ensuring that the FDR of $\empBH$ is {equal to} $\alpha/2$ for $n = m/\alpha-1$.} 
On the other hand, the variance in the FDR (blue shaded area) is smaller at fixed $n$ when $m$ increases, because the larger sample size tends to stabilize the result.

Turning to the power results, the plots show that the power of $\empBH$ also tends to that of $\BH^*$ in this sparsity regime, with also faster convergence for smaller values of $m$ (at fixed $n$). 
{This is well expected from the ``rule of thumb'' delineated in Section~\ref{sec:morefavorable} and ensuring that the transition occurs for $n\approx m/( \max(1,k) \alpha)$ where $k$ is a lower bound on the typical true discovery number of the oracle. Given the displayed results, the value of $k$ could be chosen around $(2,4,20,40)$, so that this rule would predict a transition for $n$ occurring around $(10,5,10,5)$ (top-left, top-right, bottom-left, bottom-right). Strikingly enough, the transitions indeed occur at these points in the different TDR curves. 

The sparse case where $m_1=1$ is considered in Figure~\ref{fig:puis2}, with $\mu =1$ (left column) and $\mu =3$ (right column) and a slightly increased range for $n$. Here, the oracle FDR is $0.45$ for $m=10$ and $0.495$ for $m=100$.
The observations made regarding the FDR and TDR in Figure~\ref{fig:puis1} are qualitatively the same. Moreover, in the sparse case, the convergence to the asymptotic regime is slower than in the dense case, while increasing $m$ for fixed $n$ slows down more significantly the convergence  than in the dense case. 
{This is coherent with the rule of thumb $n\approx m/( \max(1,k) \alpha)$, predicting that the transition $n$ occurs around $m/\alpha=2m$ here (only one alternative here). In addition, it is apparent on the plots that the value of the transition $n$ predicted  by this rule turns out to be particularly well adjusted, at least in this simulation setup.}

Finally, Figure~\ref{fig:puis3} compares the FDR and TDR of the procedures $\BH^*$ and $\empBH$ for larger values of $m$ and $n$  and $\alpha = 0.2$. We fix $m=10^3$ and the size of the NTS ranges from $n=1$ to $5\times 10^4$. In each plot, we see that the performances of $\empBH$ indeed increase with $n$. Despite the increased signal  amplitude in the sparse case, the situation is more difficult both in terms of convergence (which is slower) and of variance in the FDP and TDP (which are larger).
{Interestingly, this corroborates again the rule of thumb predicting a transition $n$ around $20$ and $625$ (for the choices $k\approx 250$ and $k\approx 8$) for the dense and sparse situations, respectively.}

\subsection{Additional experiments}
{
Section~\ref{add:num} presents the following additional experiments: 
first, Section~\ref{sec:naive} presents a comparison with the naive procedures $\empBY$ and $\empBHsplit$. There are both shown to be over-conservative and much less powerful than $\empBH$.
Second, {a case study with a Student distribution, leading to similar conclusions, is presented in Section~\ref{Sec:nonG}.}
Third, Section~\ref{sec:numsmalln} is devoted to simulations for very small values of $n$ ($n=5$ or $10$) with increasing values of $m$:  it shows that  $\empBH$ can achieve oracle performances in that dense case, regardless of $m$.
}

\section{{Application}}\label{sec:appli}

 One of the major scientific goals of the MUSE integral field spectrograph, which is installed on one of the 8 m telescopes at the Very Large Telescope in Chile, is the detection of distant and consequently ultra faint galaxies in the early Universe. MUSE delivers 3-dimensional datacubes (two spatial dimensions and one spectral dimension) composed of images taken in different wavelengths channels of the visible spectrum. The  values of the data samples correspond to light fluxes. Ordinary datacubes are composed  with a pile of $300\times300$  pixels  images  in $3700$ consecutive visible  wavelengths, leading to more than $300$ millions voxels. 
 
After multiple calibration and preprocessing stages, the problem of  detecting  faint galaxies boils down to a  typical needle in a haystack problem. The haystack is the datacube, which can  be considered as a discrete-valued 3-dimensional random process. This process is generated by various noise sources and by the residual perturbations of numerous bright sources. Consequently, the statistics of the random process are poorly constrained. In this haystack, each needle (there are hundreds  of them) is  a small group of connected voxels, centered on the galaxy's  position, in which the  flux locally increases. 

 A dedicated detection strategy, proposed by \cite{Origin2020} and further exploited by \cite{Cosmic2021}, consists in considering as final test statistics the 3-dimensional local maxima of the processed datacube. In the resulting testing problem, there is one null hypothesis linked to each of the $m$ local maxima, with $m$   typically in the range $[10^5,\; 10^6]$.  If we denote by $ x, y, z $ the position of a particular local maximum, we test $ H_ {0, x, y, z} $: ``There is no galaxy centred at position $ (x , y, z) $'', against $ H_ {1, x, y, z} $: ``There is one galaxy centred at this position'' and the considered  error criterion  is the FDR. 
 
As evoked above, the  distribution of the local maxima under the null hypothesis is fairly unknown. To circumvent this difficulty,  \cite {Origin2020} proposed to use the population of the opposite values of the  local minima (say, $Y_i$, in number $n$) as an independent ``proxy'' (a NTS)  for the local maxima (say, $X_i$,  in number $m$). They reported  numerical simulations  suggesting that a procedure close to the Benjamini-Hochberg procedure using $p$-values computed from this NTS controls the FDR.  This astrophysical application involves a common but unknown distribution $P_0$ under the null hypothesis and the possibility of using a NTS to improve the control of the FDR: this is clearly the setting described in \ref{sec:setting}  and in fact this application has inspired the present study. It is thus interesting to see which light the present study sheds on this initial approach. 

The sample sizes considered here are $n=2.3\times 10^6$ and $m=3.3\times 10^6$ so $n<m$ and both are large. The empirical distribution of the values of the NTS and of the test sample are shown in Figure~\ref{fig:muse}, left panel. The similarity of the two distributions in the left and central parts suggests that the NTS (blue) can serve as a useful proxy for the test sample (red). The right tail of the test sample is logically heavier owing to the presence of galaxies, which tend to shift the values of the local maxima upwards.

While the procedure proposed by \cite{Origin2020} is very close to  the $\empBH$ procedure with Algorithm \ref{algo:empBH}, it differs in the following point.
Instead of using
FDP $=\frac{{V+1}}{K}\frac{m}{n+1}$ (see step 4 of Algorithm \ref{algo:empBH}),  \cite{Origin2020} use  FDP $=\frac{{V}}{K}\frac{N_1}{N_0}$, where $N_1$ (resp. $N_0$) are  the number of voxels of the region where the local extrema  of the test sample (resp., the NTS) are computed.  Because $n$ is  large, $V$ is large as well and using $V$ instead of $V+1$ has no numerical impact in this regime.   The normalization factors are in fact very similar as well, with $\frac{m}{n+1}\approx 1.440$ and  $\frac{N_1}{N_0}\approx 1.442$. 
In effect,  it turns out that there is no numerical difference in running these two versions of Algorithm~\ref{algo:empBH}: in both cases, the procedure rejects exactly $105$ local maxima at target FDR $\alpha = 0.2$, a situation shown in the right panel of Figure~\ref{fig:muse}. 
  {The rejected local maxima in \cite{Origin2020} being the same as those rejected by $\empBH$, the discovery set inherits the properties delineated in the present work: first, the FDR control is established from Theorem~\ref{thFDR}, because the only assumption made (Assumption \eqref{indep}) is likely to hold due to the important dimension reduction made in the dataset when focusing on the local maxima/minima. Second, we have $n\approx m$ while both are large. This means that we are just at the border of the boundary identified in Section~\ref{sec:boundary}, so the theory is silent for this case. Nevertheless, the distribution of the data exhibits some minimum amount of signal, perhaps  $k\gtrsim 50$ fairly detectable alternatives. Hence, the refined upper-bound given in Proposition~\ref{extension} can also be applied: since the training-to-test ratio $n/m$ is above the boundary, that is, $n/m\approx 1$ is much larger than $1/(k \alpha)\lesssim 0.1$, the power of $\empBH$ should be close to the one of the oracle for this data set. }

To conclude, the present paper illustrates that $\empBH$, together with our theoretical findings, delivers interpretable and useful  results for common practice.   
Meanwhile, it validates the use of the  procedure proposed
in \cite{Origin2020} on this particular data set. 

\begin{figure}
    \centering
    \subfigure[]{\includegraphics[width=.49\textwidth]{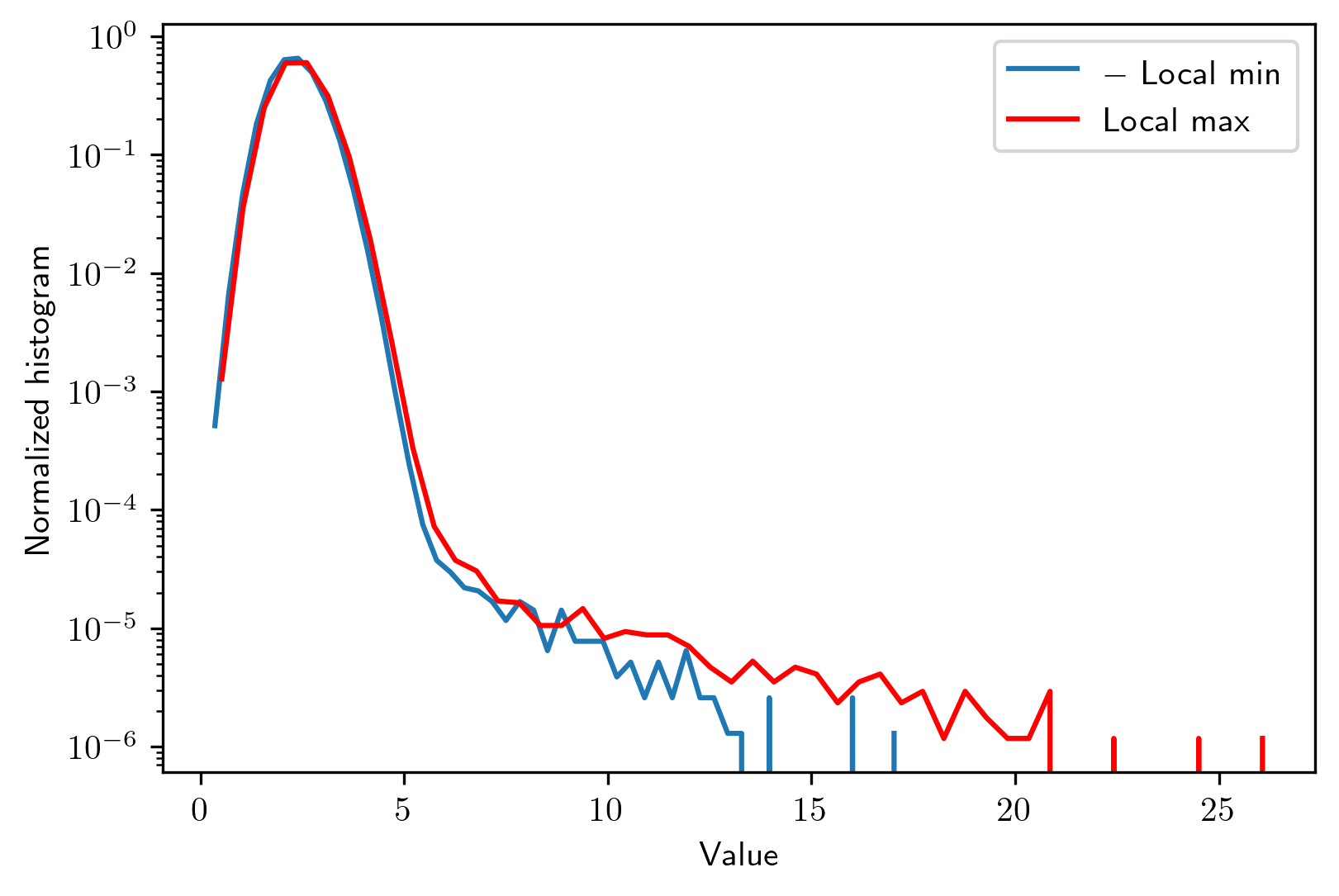}} 
  \subfigure[]{\includegraphics[width=0.49\textwidth]{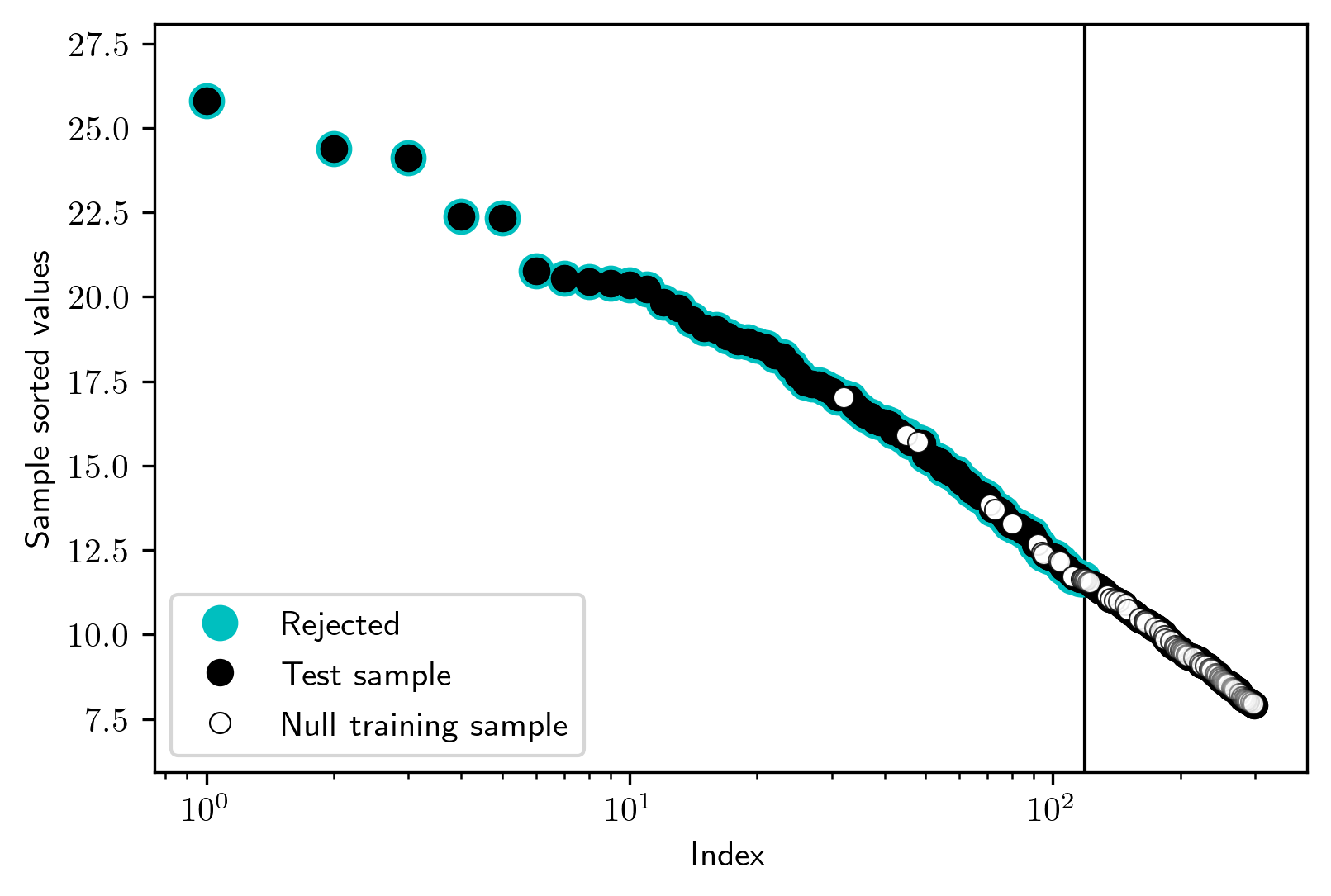}} 
    \caption{\label{fig:muse} MUSE application example. (a) Empirical distributions   of the values of the NTS (the $Y_i$, computed as the opposite of the local minima, in blue) and of the test sample (the $X_i$, local maxima, in red). 
     (b) { Results of Algorithm 1 (same color code as in Figure~{\ref{fig:algo}).     
    Only the $700$ largest sample values are shown. The black  vertical line indicates the rejection  {threshold $K=105$} and  $14$ samples  of the NTS are above this threshold.}} }   \label{fig:foobar}
\end{figure}

\section{Conclusion and discussion}\label{sec:conclusion}

\subsection{Summary}

In a nutshell, this paper evaluated how classical multiple testing methodology can generalize when replacing the knowledge of the null distribution $P_0$ by examples $Y_1,\dots,Y_n$ following this null. 
While this situation is very frequent in practice, it has only been scarcely studied so far and this paper contributed to fill this gap.
The FDR control guarantee holds whatever $n,m$, with no assumption on $P_0$ and for any marginal alternative, with a bound $\alpha m_0/m$ (achieved when $\alpha(n+1)/m$ is an integer), which is similar to the result obtained in the original work of \cite{BH1995} in case where $P_0$ is known. In addition, the power is comparable to the one of the oracle when $n\gtrsim  m/(\alpha\max(1,k))$, where $k$ is a confidence lower bound on the number of true discoveries made by the oracle. This ``rule of thumb'' has been both validated by theory and numerical experiments. 
Finally, we demonstrated that our work brought a theoretical support and thus more interpretability in a  worked-out application to recent breakthrough findings in astrophysics. 
{In practice, our ``rule of thumb'' $n\gtrsim  m/(\alpha\max(1,k))$ can be used as follows: if the user has no strong prior belief in a minimum number $k$ of discoveries, choosing $k=0$ might be safer, which leads to the condition $n\gtrsim  m/\alpha$. By contrast, if $k$ can be accurately guessed a priori, the less demanding condition $n\gtrsim  m/\alpha k$ can be opted for.}

This work also completed the picture by exhibiting a theoretical intrinsic limitation of the semi-supervised multiple testing setting when the null training sample is not populated enough. It is impossible to control the FDR while mimicking the oracle power for $n\lesssim m$ when letting the sparsity and the distribution of the alternative arbitrary. This delineates a setting-intrinsic phase transition at $n\asymp m$. 

\subsection{Optimality of $\BH^*$}\label{sec:BHoptimal}

In the past literature, a numerous of works proposed approaches that improve, sometimes substantially, the baseline BH procedure (as local FDR methods listed in introduction). Hence, a common belief is that the BH procedure is well known to be conservative and suboptimal when controlling the FDR. This belief makes the aim of mimicking the performance of the oracle BH procedure (as in Sections~\ref{powerresult}~and~\ref{sec:optimality}) somewhat questionable. However, we argue that this belief is not justified  when the test statistic used before applying BH algorithm is suitably chosen, typically using a likelihood ratio or a local FDR transformation. This is shown in particular with simulations in the setting of Appendix~\ref{sec:MC}, where $\BH^*$ improves over a local FDR method, itself well known to enjoy optimality properties \citep{CSWW2019}. In a nutshell, the possible conservativeness of BH procedure when controlling the FDR is not due to the BH algorithm {\it per se} but rather to the test statistic used as entries of this algorithm. To come back to our framework, the test statistic is assumed to be fixed once for all in our work. 
Hence, {\it given the chosen test statistic}, $\BH^*$ is close to be optimal and the aim considered in  Sections~\ref{powerresult}~and~\ref{sec:optimality} perfectly makes sense.



\subsection{Future work}

Given that semi-supervised multiple testing setting is versatile, our work raises a number of new perspectives. For instance, in recent machine learning, this setting conveniently bypasses model assumptions on $P_0$ and only needs a number of null examples, that can be generated by a suitable ``blackbox''. Nevertheless, in order to avoid potential bias in the null training sample, this blackbox should be properly calibrated with significant prior calibrations and preprocessing steps. While building such an approach deserves an entire devoted study, we anticipate that studying the robustness of the procedure  $\empBH$ with respect the NTS is a key point: what about the case where $Y_1,\dots,Y_n$ are i.i.d. $\sim P_0'$ with $P_0'\approx P_0$?

Another avenue for future work is to decline recent advances in multiple testing into this semi-supervised setting. 
For instance, while $\empBH$ is devoted to the FDR criterion, an interesting and challenging issue is to design semi-supervised counterparts suitable for other criteria, as FDX \cite{GW2004}, online FDR \cite{FS2008,xu2021dynamic} or post hoc bounds \cite{GW2006,GS2011}. In particular, since the variability of the FDP of $\empBH$ is increased by the NTS, considering criteria accounting for this effect seems particularly interesting. 
Since various dependence assumptions are used in such studies, we also expect that our main assumption \eqref{exch} can be relaxed in some of these frameworks.

Finally, proper calibrations of the individual tests sometimes require to consider hypothesis-dependent null distributions, that is, null distributions $P_{0,i}$ that depend on $i\in\{1,\dots,m\}$ (see, e.g., \citealp{SMB2017,Sulis2020} for a concrete example). Since $m$ null training samples should be considered in that case, it poses a complexity issue and generalizing our result to this setting is both 
theoretically challenging and useful to support or improve procedures used in common practice.

\section*{Acknowledgements}
This work has been supported by ANR-16-CE40-0019 (SansSouci), ANR-17-CE40-0001 (BASICS) and by the GDR ISIS through the ``projets exploratoires" program (project TASTY).
We are grateful to Sabine Houssaye for her help when proving Lemma~\ref{proof:algebre} and to Guillaume Lecu\'e for helpful comments. 

\bibliographystyle{apalike}
\bibliography{biblio}

\begin{thebibliography}{}

\bibitem[Abraham et~al., 2021]{abraham2021multiple}
Abraham, K., Castillo, I., and Gassiat, E. (2021).
\newblock Multiple testing in nonparametric hidden markov models: An empirical
  bayes approach.
\newblock {\em arXiv preprint arXiv:2101.03838}.

\bibitem[Arlot et~al., 2010]{ABR2010a}
Arlot, S., Blanchard, G., and Roquain, E. (2010).
\newblock Some nonasymptotic results on resampling in high dimension. {I}.
  {C}onfidence regions.
\newblock {\em Ann. Statist.}, 38(1):51--82.

\bibitem[Azriel and Schwartzman, 2015]{AS2015}
Azriel, D. and Schwartzman, A. (2015).
\newblock The empirical distribution of a large number of correlated normal
  variables.
\newblock {\em Journal of the American Statistical Association},
  110(511):1217--1228.

\bibitem[Bacon et~al., 2021]{Cosmic2021}
Bacon, R., {Mary, D.}, {Garel, T.}, {Blaizot, J.}, {Maseda, M.}, {Schaye, J.},
  {Wisotzki, L.}, {Conseil, S.}, {Brinchmann, J.}, {Leclercq, F.},
  {Abril-Melgarejo, V.}, {Boogaard, L.}, {Bouch\'e, N. F.}, {Contini, T.},
  {Feltre, A.}, {Guiderdoni, B.}, {Herenz, C.}, {Kollatschny, W.}, {Kusakabe,
  H.}, {Matthee, J.}, {Michel-Dansac, L.}, {Nanayakkara, T.}, {Richard, J.},
  {Roth, M.}, {Schmidt, K. B.}, {Steinmetz, M.}, {Tresse, L.}, {Urrutia, T.},
  {Verhamme, A.}, {Weilbacher, P. M.}, {Zabl, J.}, and {Zoutendijk, S. L.}
  (2021).
\newblock The muse extremely deep field: The cosmic web in emission at high
  redshift.
\newblock {\em A\&A}, 647:A107.

\bibitem[Barber and Cand\`es, 2015]{BC2015}
Barber, R.~F. and Cand\`es, E.~J. (2015).
\newblock Controlling the false discovery rate via knockoffs.
\newblock {\em Ann. Statist.}, 43(5):2055--2085.

\bibitem[{Barber} and {Cand\`es}, 2019]{barber2019knockoff}
{Barber}, R.~F. and {Cand\`es}, E.~J. (2019).
\newblock {A knockoff filter for high-dimensional selective inference}.
\newblock {\em {Ann. Stat.}}, 47(5):2504--2537.

\bibitem[{Barber} et~al., 2020]{Barber2020robust}
{Barber}, R.~F., {Cand\`es}, E.~J., and {Samworth}, R.~J. (2020).
\newblock {Robust inference with knockoffs}.
\newblock {\em {Ann. Stat.}}, 48(3):1409--1431.

\bibitem[Bates et~al., 2020]{bates2020metropolized}
Bates, S., Cand{\`e}s, E., Janson, L., and Wang, W. (2020).
\newblock Metropolized knockoff sampling.
\newblock {\em Journal of the American Statistical Association}, pages 1--15.

\bibitem[Bates et~al., 2021]{bates2021testing}
Bates, S., Cand\`es, E., Lei, L., Romano, Y., and Sesia, M. (2021).
\newblock Testing for outliers with conformal p-values.

\bibitem[Benjamini and Hochberg, 1995]{BH1995}
Benjamini, Y. and Hochberg, Y. (1995).
\newblock Controlling the false discovery rate: a practical and powerful
  approach to multiple testing.
\newblock {\em J. Roy. Statist. Soc. Ser. B}, 57(1):289--300.

\bibitem[Benjamini and Yekutieli, 2001]{BY2001}
Benjamini, Y. and Yekutieli, D. (2001).
\newblock The control of the false discovery rate in multiple testing under
  dependency.
\newblock {\em Ann. Statist.}, 29(4):1165--1188.

\bibitem[Besag and Clifford, 1991]{besag1991sequential}
Besag, J. and Clifford, P. (1991).
\newblock Sequential monte carlo p-values.
\newblock {\em Biometrika}, 78(2):301--304.

\bibitem[Blanchard et~al., 2010]{BLS2010}
Blanchard, G., Lee, G., and Scott, C. (2010).
\newblock Semi-supervised novelty detection.
\newblock {\em J. Mach. Learn. Res.}, 11:2973--3009.

\bibitem[Blanchard and Roquain, 2008]{BR2008}
Blanchard, G. and Roquain, E. (2008).
\newblock Two simple sufficient conditions for {FDR} control.
\newblock {\em Electron. J. Stat.}, 2:963--992.

\bibitem[Cai and Sun, 2009]{CS2009}
Cai, T.~T. and Sun, W. (2009).
\newblock Simultaneous testing of grouped hypotheses: finding needles in
  multiple haystacks.
\newblock {\em J. Amer. Statist. Assoc.}, 104(488):1467--1481.

\bibitem[Cai et~al., 2019]{CSWW2019}
Cai, T.~T., Sun, W., and Wang, W. (2019).
\newblock Covariate-assisted ranking and screening for large-scale two-sample
  inference.
\newblock In {\em Royal Statistical Society}, volume~81.

\bibitem[{Cand\`es} et~al., 2018]{candes2018panning}
{Cand\`es}, E., {Fan}, Y., {Janson}, L., and {Lv}, J. (2018).
\newblock {Panning for gold: `model-X' knockoffs for high dimensional
  controlled variable selection}.
\newblock {\em {J. R. Stat. Soc., Ser. B, Stat. Methodol.}}, 80(3):551--577.

\bibitem[Carpentier et~al., 2021]{CDRV2021}
Carpentier, A., Delattre, S., Roquain, E., and Verzelen, N. (2021).
\newblock Estimating minimum effect with outlier selection.
\newblock {\em Annals of Statistics}, 49(1):272--294.

\bibitem[Choquet et~al., 2018]{Choquet_2018}
Choquet, {\'{E}}., Bryden, G., Perrin, M.~D., Soummer, R., Augereau, J.-C.,
  Chen, C.~H., Debes, J.~H., Gofas-Salas, E., Hagan, J.~B., Hines, D.~C.,
  Mawet, D., Morales, F., Pueyo, L., Rajan, A., Ren, B., Schneider, G., Stark,
  C.~C., and Wolff, S. (2018).
\newblock {HD} 104860 and {HD} 192758: Two debris disks newly imaged in
  scattered light with {theHubble} space telescope.
\newblock {\em The Astrophysical Journal}, 854(1):53.

\bibitem[Davison and Hinkley, 1997]{davison1997bootstrap}
Davison, A.~C. and Hinkley, D.~V. (1997).
\newblock {\em Bootstrap methods and their application}.
\newblock Number~1. Cambridge university press.

\bibitem[Dunnett, 1955]{dunnett1955multiple}
Dunnett, C.~W. (1955).
\newblock A multiple comparison procedure for comparing several treatments with
  a control.
\newblock {\em Journal of the American Statistical Association},
  50(272):1096--1121.

\bibitem[{Efron}, 2004]{Efron2004}
{Efron}, B. (2004).
\newblock {Large-scale simultaneous hypothesis testing: the choice of a null
  hypothesis.}
\newblock {\em {J. Am. Stat. Assoc.}}, 99(465):96--104.

\bibitem[Efron, 2007]{Efron2007b}
Efron, B. (2007).
\newblock Doing thousands of hypothesis tests at the same time.
\newblock {\em Metron - International Journal of Statistics}, LXV(1):3--21.

\bibitem[Efron, 2008]{Efron2008}
Efron, B. (2008).
\newblock Microarrays, empirical {B}ayes and the two-groups model.
\newblock {\em Statist. Sci.}, 23(1):1--22.

\bibitem[{Efron}, 2009]{Efron2009b}
{Efron}, B. (2009).
\newblock {Empirical Bayes estimates for large-scale prediction problems.}
\newblock {\em {J. Am. Stat. Assoc.}}, 104(487):1015--1028.

\bibitem[Efron et~al., 2001]{ETST2001}
Efron, B., Tibshirani, R., Storey, J.~D., and Tusher, V. (2001).
\newblock Empirical {B}ayes analysis of a microarray experiment.
\newblock {\em J. Amer. Statist. Assoc.}, 96(456):1151--1160.

\bibitem[Finner and Strassburger, 2007]{FS2007b}
Finner, H. and Strassburger, K. (2007).
\newblock Step-up related simultaneous confidence intervals for mcc and mcb.
\newblock {\em Biometrical Journal}, 49(1):40--51.

\bibitem[Fisher, 1935]{Fish1935}
Fisher, R.~A. (1935).
\newblock {\em The Design of Experiments.}
\newblock Oliver and Boyd, Edinburgh.

\bibitem[Fithian and Lei, 2020]{fithian2020conditional}
Fithian, W. and Lei, L. (2020).
\newblock Conditional calibration for false discovery rate control under
  dependence.

\bibitem[{Foster} and {Stine}, 2008]{FS2008}
{Foster}, D.~P. and {Stine}, R.~A. (2008).
\newblock {\(\alpha\)-investing: a procedure for sequential control of expected
  false discoveries}.
\newblock {\em {J. R. Stat. Soc., Ser. B, Stat. Methodol.}}, 70(2):429--444.

\bibitem[{Gandy} and {Hahn}, 2014]{GH2014}
{Gandy}, A. and {Hahn}, G. (2014).
\newblock {MMCTest -- a safe algorithm for implementing multiple Monte Carlo
  tests}.
\newblock {\em {Scand. J. Stat.}}, 41(4):1083--1101.

\bibitem[Genovese and Wasserman, 2004]{GW2004}
Genovese, C. and Wasserman, L. (2004).
\newblock A stochastic process approach to false discovery control.
\newblock {\em Ann. Statist.}, 32(3):1035--1061.

\bibitem[Genovese and Wasserman, 2006]{GW2006}
Genovese, C.~R. and Wasserman, L. (2006).
\newblock Exceedance control of the false discovery proportion.
\newblock {\em J. Amer. Statist. Assoc.}, 101(476):1408--1417.

\bibitem[Goeman and Solari, 2011]{GS2011}
Goeman, J.~J. and Solari, A. (2011).
\newblock Multiple testing for exploratory research.
\newblock {\em Statist. Sci.}, 26(4):584--597.

\bibitem[Goodfellow et~al., 2014]{Good2014}
Goodfellow, I., Pouget-Abadie, J., Mirza, M., Xu, B., Warde-Farley, D., Ozair,
  S., Courville, A., and Bengio, Y. (2014).
\newblock Generative adversarial nets.
\newblock In Ghahramani, Z., Welling, M., Cortes, C., Lawrence, N., and
  Weinberger, K.~Q., editors, {\em Advances in Neural Information Processing
  Systems}, volume~27. Curran Associates, Inc.

\bibitem[{Guo} and {Peddada}, 2008]{GP2008}
{Guo}, W. and {Peddada}, S. (2008).
\newblock {Adaptive choice of the number of bootstrap samples in large scale
  multiple testing}.
\newblock {\em {Stat. Appl. Genet. Mol. Biol.}}, 7(1):19.
\newblock Id/No 13.

\bibitem[Heller and Yekutieli, 2014]{heller2014}
Heller, R. and Yekutieli, D. (2014).
\newblock Replicability analysis for genome-wide association studies.
\newblock {\em Ann. Appl. Stat.}, 8(1):481--498.

\bibitem[Hemerik et~al., 2019]{HSG2019}
Hemerik, J., Solari, A., and Goeman, J.~J. (2019).
\newblock {Permutation-based simultaneous confidence bounds for the false
  discovery proportion}.
\newblock {\em Biometrika}, 106(3):635--649.

\bibitem[Hsu, 1996]{hsu1996multiple}
Hsu, J. (1996).
\newblock {\em Multiple comparisons: theory and methods}.
\newblock CRC Press.

\bibitem[Katsevich and Sabatti, 2019]{katsevich2019multilayer}
Katsevich, E. and Sabatti, C. (2019).
\newblock Multilayer knockoff filter: Controlled variable selection at multiple
  resolutions.
\newblock {\em The annals of applied statistics}, 13(1):1.

\bibitem[Kingma and Welling, 2014]{Kingma2014}
Kingma, D.~P. and Welling, M. (2014).
\newblock Auto-encoding variational bayes.
\newblock In Bengio, Y. and LeCun, Y., editors, {\em 2nd International
  Conference on Learning Representations, {ICLR} 2014, Banff, AB, Canada, April
  14-16, 2014, Conference Track Proceedings}.

\bibitem[Lin, 2005]{lin2005efficient}
Lin, D. (2005).
\newblock An efficient monte carlo approach to assessing statistical
  significance in genomic studies.
\newblock {\em Bioinformatics}, 21(6):781--787.

\bibitem[Liu and Zheng, 2018]{liu2018auto}
Liu, Y. and Zheng, C. (2018).
\newblock Auto-encoding knockoff generator for fdr controlled variable
  selection.
\newblock {\em arXiv preprint arXiv:1809.10765}.

\bibitem[Mary et~al., 2020]{Origin2020}
Mary, D., Bacon, R., Conseil, S., Piqueras, L., and Schutz, A. (2020).
\newblock {ORIGIN}: {B}lind detection of faint emission line galaxies in muse
  datacubes.
\newblock {\em A\&A}, 635:A194.

\bibitem[Nguyen et~al., 2020]{nguyen2020aggregation}
Nguyen, T.-B., Chevalier, J.-A., Thirion, B., and Arlot, S. (2020).
\newblock Aggregation of multiple knockoffs.
\newblock In {\em International Conference on Machine Learning}, pages
  7283--7293. PMLR.

\bibitem[Padilla and Bickel, 2012]{PB2012}
Padilla, M. and Bickel, D.~R. (2012).
\newblock Estimators of the local false discovery rate designed for small
  numbers of tests.
\newblock {\em Stat. Appl. Genet. Mol. Biol.}, 11(5):Art. 4, front matter+39.

\bibitem[Phipson and Smyth, 2010]{phipson2010permutation}
Phipson, B. and Smyth, G.~K. (2010).
\newblock Permutation p-values should never be zero: calculating exact p-values
  when permutations are randomly drawn.
\newblock {\em Statistical applications in genetics and molecular biology},
  9(1).

\bibitem[Romano and Wolf, 2005]{RW2005}
Romano, J.~P. and Wolf, M. (2005).
\newblock Exact and approximate stepdown methods for multiple hypothesis
  testing.
\newblock {\em J. Amer. Statist. Assoc.}, 100(469):94--108.

\bibitem[Romano and Wolf, 2007]{RW2007}
Romano, J.~P. and Wolf, M. (2007).
\newblock Control of generalized error rates in multiple testing.
\newblock {\em Ann. Statist.}, 35(4):1378--1408.

\bibitem[Roquain and Verzelen, 2020a]{RVvignette2020}
Roquain, E. and Verzelen, N. (2020a).
\newblock False discovery rate control with unknown null distribution:
  illustrations on real data sets.
\newblock https://github.com/eroquain/empiricalnull/blob/main/vignette.pdf.

\bibitem[Roquain and Verzelen, 2020b]{roquain2020false}
Roquain, E. and Verzelen, N. (2020b).
\newblock False discovery rate control with unknown null distribution: is it
  possible to mimic the oracle?

\bibitem[Sandve et~al., 2011]{sandve2011sequential}
Sandve, G.~K., Ferkingstad, E., and Nyg{\aa}rd, S. (2011).
\newblock Sequential monte carlo multiple testing.
\newblock {\em Bioinformatics}, 27(23):3235--3241.

\bibitem[Sarkar and Tang, 2021]{sarkar2021adjusting}
Sarkar, S.~K. and Tang, C.~Y. (2021).
\newblock Adjusting the benjamini-hochberg method for controlling the false
  discovery rate in knockoff assisted variable selection.
\newblock {\em arXiv preprint arXiv:2102.09080}.

\bibitem[Schwartzman, 2010]{Sch2010}
Schwartzman, A. (2010).
\newblock Comment: ``{C}orrelated {$z$}-values and the accuracy of large-scale
  statistical estimates''.
\newblock {\em J. Amer. Statist. Assoc.}, 105(491):1059--1063.

\bibitem[Stephens, 2017]{Ste2017}
Stephens, M. (2017).
\newblock False discovery rates: a new deal.
\newblock {\em Biostatistics}, 18(2):275--294.

\bibitem[Sulis et~al., 2017]{SMB2017}
Sulis, S., Mary, D., and Bigot, L. (2017).
\newblock A study of periodograms standardized using training datasets and
  application to exoplanet detection.
\newblock {\em IEEE Transactions on Signal Processing}, 65(8):2136--2150.

\bibitem[Sulis et~al., 2020]{Sulis2020}
Sulis, S., Mary, D., and Bigot, L. (2020).
\newblock 3{D} magneto-hydrodynamical simulations of stellar convective noise
  for improved exoplanet detection - {I}. {C}ase of regularly sampled radial
  velocity observations.
\newblock {\em A\&A}, 635:A146.

\bibitem[Sun and Stephens, 2018]{SS2018}
Sun, L. and Stephens, M. (2018).
\newblock Solving the empirical bayes normal means problem with correlated
  noise.

\bibitem[{Sun} and {Cai}, 2007]{SC2007}
{Sun}, W. and {Cai}, T.~T. (2007).
\newblock {Oracle and adaptive compound decision rules for false discovery rate
  control}.
\newblock {\em {J. Am. Stat. Assoc.}}, 102(479):901--912.

\bibitem[Sun and Cai, 2009]{SC2009}
Sun, W. and Cai, T.~T. (2009).
\newblock Large-scale multiple testing under dependence.
\newblock {\em J. R. Stat. Soc. Ser. B Stat. Methodol.}, 71(2):393--424.

\bibitem[Weinstein et~al., 2017]{weinstein2017power}
Weinstein, A., Barber, R., and Cand{\`e}s, E. (2017).
\newblock A power and prediction analysis for knockoffs with lasso statistics.

\bibitem[Westfall and Young, 1993]{WY1993}
Westfall, P.~H. and Young, S.~S. (1993).
\newblock {\em Resampling-Based Multiple Testing}.
\newblock Wiley.
\newblock Examples and Methods for $P$- Value Adjustment.

\bibitem[Xu and Ramdas, 2021]{xu2021dynamic}
Xu, Z. and Ramdas, A. (2021).
\newblock Dynamic algorithms for online multiple testing.

\bibitem[Yang et~al., 2021]{yang2021bonus}
Yang, C.-Y., Lei, L., Ho, N., and Fithian, W. (2021).
\newblock Bonus: Multiple multivariate testing with a data-adaptive test
  statistic.

\bibitem[Zhang et~al., 2019]{zhang_adaptive_2019}
Zhang, M.~J., Zou, J., and Tse, D. (2019).
\newblock Adaptive {Monte} {Carlo} {Multiple} {Testing} via {Multi}-{Armed}
  {Bandits}.
\newblock {\em arXiv:1902.00197 [cs, math, q-bio, stat]}.
\newblock arXiv: 1902.00197.

\end{thebibliography}
\appendix

\section{By-product 1: Blackbox BH procedure}\label{sec:MC}

\subsection{Setting and procedure}

In this section, we consider the same formal setting and notation as in Section~\ref{sec:setting}, except that the test statistics $X_1,\dots,X_m$ are given along with a ``blackbox sampler'' able to produce i.i.d. realizations of the null $P_0$, even if $P_0$ is not known. As in the motivations described in Section~\ref{sec:background}, such a blackbox can come from an external code implemented by an expert of the application domain, or from a machine learning program that has been sufficiently trained. 
Our work easily allows to design a multiple testing inference in that situation. 
  Namely, Algorithm~\ref{algo:MCBH} below can be used to produce a sampled BH procedure, that we call the Blackbox BH procedure (bbBH). 
By Theorem~\ref{thFDR}, the bbBH procedure achieves an FDR equal to $\alpha m_0/m$ (when $\alpha$ is a rational number), provided that $(X_i,i\in\cH_0)$ are i.i.d. $\sim P_0$ and independent of  $(X_i,i\notin\cH_0)$. 
Also, since $n$ is chosen so that $(n+1)\alpha/m\geq 1$, {it is just above the boundary put forward in Section~\ref{powerresult}, which might indicate that the power of bbBH should be comparable to that of the oracle. }

\begin{algorithm}
\KwData{$X=(X_1,\dots,X_m)\in \R^{m}$, a nominal level $\alpha\in (0,1)$ (assumed to be a rational number) and a blackbox sampler of the null distribution $P_0$}
 \begin{enumerate}
 \item Choose $n\geq 1$ the smallest integer such that $(n+1)\alpha/m$ is an integer
 \item Sample $(Y_1,\dots,Y_n)$ i.i.d. according to the null distribution $P_0$
\item Apply the semi-supervised BH procedure $\empBH_\alpha$ to $Z=(Y,X)$, see Algorithm~\ref{algo:empBH}
\end{enumerate}
\KwResult{Reject the nulls in the set $\empBH_\alpha$.}
\caption{Blackbox BH procedure }\label{algo:MCBH}
\end{algorithm}

\subsection{Illustration with simultaneous likelihood ratio tests}\label{sec:LRT}

To illustrate further the interest of the bbBH procedure, we consider in this section the problem of controlling the FDR while choosing the best individual test statistics.
To this end, let us consider the common setup where we   
observe $m$ independent measurements $T_1,\dots,T_m\in \R$, with $T_i$ either distributed as a null distribution $G_0$ or as an alternative distribution $G_1$, where $G_0$ and $G_1$ are {\it known} distribution with densities $g_0$ and $g_1$, respectively. For each $i\in \{1,\dots,m\}$,  
we consider a likelihood ratio test of the null hypothesis $H_{0,i}$: ``$T_i\sim G_0$'' against the alternative $H_{1,i}$: ``$T_i\sim G_1$''. It rejects the null whenever $\{X_i\geq c\}$ for $X_i=g_1(T_i)/g_0(T_i)$ (with the convention $X_i=+\infty$ if $g_0(T_i)=0$) and some constant $c>0$ such that $\bar{F}_0(c)=\alpha$ with 
\begin{equation}\label{F0barbb}
\bar{F}_0(t)=\int_\R \ind{g_1(u)>t g_0(u)}g_0(u) du, \:\: t\geq 0.
\end{equation}
 Denote $P_0$ the distribution of $X_i$ under $G_0$ and assume that $\bar{F}_0$ is continuous and decreasing on the support of $P_0$.  
{The oracle BH procedure $\BH^*$, which is not accessible in general, can be nevertheless approximated in this setting via a numerical 
approximation of the function $\bar{F}_0$. 
By contrast, we can also build a ``blackbox'' that generates realizations of $P_0$ by simulating $T_1,\dots,T_n$ i.i.d. $\sim G_0$ and then letting $Y_i=g_1(T_i)/g_0(T_i)$, $1\leq i\leq n$. }
Hence, we can apply  the bbBH procedure (Algorithm~\ref{algo:MCBH}) to control the FDR 
at level $\alpha$ in this model (and even have an FDR equal to $\alpha m_0/m$), while having a power close to the one of the oracle.
 
 For comparisons, we also introduce two other procedures: first, the BH procedure directly applied to the original test statistics $T_i$'s (with respect to the known null $G_0$), which is referred to as BH$0$ below. 
Compared to bbBH, BH$0$ has the advantage to be not-random. 
However, since the individual tests based on the $T_i$ are less powerful than those based on the likelihood ratio $X_i$, bbBH is in general more powerful than BH$0$. 
The second procedure is the classical FDR controlling method based on the local FDR values \cite{ETST2001,SC2007}, denoted by ``locfdr'', which can be used specifically for this example, see Section~\ref{sec:numbbBH} for more details.

The performances of these procedure are illustrated by a numerical experiments in Section~\ref{sec:numbbBH}. The conclusions of this experiment are as follows:
\begin{itemize}
\item All procedures correctly control the FDR; 
\item As expected bbBH, $\BH^*$ and locfdr have better power than BH$0$; 
\item The two procedures locfdr and bbBH both mimic the power of the oracle $\BH^*$, although bbBH is better adjusted for $m$ small, 
while locfdr is slightly less variant. 
\end{itemize}
Overall, this section validates the use of bbBH in a ``toy blackbox setting'' where alternative procedures can be employed. This suggests that bbBH will perform favorably in general blackbox settings for which no such alternative exists.

\section{By-product 2: the randomized BH procedure}\label{sec:randBH}

Let us consider, only for the present section, the usual framework where the null distribution is known and no NTS are given. In particular, $\BH^*$ boils down to the usual BH procedure.

Recall that an important part of multiple testing literature is devoted to find procedures that control rigorously the FDR at level $\alpha$ while maximizing the power. We emphasize that, in this framework, having an FDR equal to $\alpha+10^{-10}$ (say) is not allowed: the inequality $\FDR\leq \alpha$ must hold under any configuration of the model, which is particularly challenging when negative dependences are possible. For instance, we refer to the very recent work of \cite{fithian2020conditional} (see also references therein), that aims at modifying the BH procedure in order to control the FDR under negative dependence.
The point of this section is to point out that Theorem~\ref{thFDR} and Example~\ref{sec:robutnegative} allow to solve this problem in a simple way for some (admittedly specific) dependence structure.  

Assume that $X$ is an $m$-dimensional Gaussian equi-correlated vector with individual variances equal to  $1$ and with {\it known} covariance $\rho\in [-1/m,0)$. Consider $n=n(\rho,m)\geq 1$ the largest integer so that $\rho\geq -1/(n+m-1)$, that is, $n=\lfloor -\rho^{-1}-m+1\rfloor$ and generate a $n$-sample $Y_1,\dots,Y_n$ such that $(Y,X)$ is Gaussian equi-correlated $\rho$. This can be done easily via Proposition~\ref{prop:equialgo}. Then Theorem~\ref{thFDR} provides that the procedure $\empBH_\alpha$ controls the FDR at level $\alpha$. Here, since the NTS is generated by the user, this procedure can be seen as a {\it randomized BH procedure}, (randBH in short). Algorithm~\ref{algo:randBH} gives the full steps to implement randBH. 
In addition, our rule of thumb suggests that RandBH has a power comparable to that of BH when $\rho\gtrsim -1/(m+m/(\alpha \max(1,k)))$, that is, for configurations close enough to the independent case.

\begin{algorithm}
\KwData{$X=(X_1,\dots,X_m)\in \R^{m}$, $\rho\in [-1/m,0)$, $\alpha$}
 \begin{enumerate}
 \item Compute $n\geq 1$ the largest integer so that $\rho\geq -1/(n+m-1)$, that is, $n=\lfloor -\rho^{-1}-m+1\rfloor$
\item Let $T=X$
\item For $k$ from $m$ to $n+m-1$:
\begin{itemize}
\item draw $U\sim \mathcal{N}(0,1)$ independently of the rest
\item 
let $
T_{k+1} = \frac{\rho}{1+(k-1)\rho} (T_1+\dots+T_k) + \left(1-k\frac{\rho^2}{1+(k-1)\rho}  \right)^{1/2} U
$
\item let $T=(T_1,\dots,T_{k+1})$
\end{itemize}
\item Let $Y=(T_{m+1},\dots,T_{n+m})$
\item Apply the semi-supervised BH procedure $\empBH_\alpha$ to $Z=(Y,X)$, see Algorithm~\ref{algo:empBH}
\end{enumerate}
\KwResult{Reject the nulls in the set $\empBH_\alpha$.}
\caption{Randomized BH procedure}\label{algo:randBH}
\end{algorithm}

Also, we would like to make a disclaimer: {we do not pretend that RandBH is applicable in general practice, at least under the current form, because it is linked to a too specific dependence structure}. Rather, the message is that randomization (plus using $p$-values biased upwards) can help the BH procedure to be more robust with respect to negative dependencies. We think that this intriguing side result is an important proof of concept.

Finally, this phenomenon can be derived for other negative dependence structures: however, the reachable distributions of $(X_i,i\in \cH_0)$ should necessarily be expressible as a marginal of a larger vector $(Y_1,\dots,Y_n,X_i,i\in\cH_0)$ that is exchangeable in order to satisfy \eqref{exch}.

\end{document}